\documentclass[a4paper]{article}

\usepackage{amsmath,amssymb,amsthm,a4wide,mathrsfs}
\usepackage{graphicx,color,caption,subcaption}
\usepackage{algorithm,algpseudocode}
\usepackage[utf8]{inputenc}
\usepackage[T1]{fontenc}
\usepackage{tikz}
\usepackage{hyperref}
\hypersetup{
    bookmarksopen=true,
    colorlinks=true,
    linkcolor=blue,
    citecolor=blue,
    allcolors = {blue},
}

\graphicspath{{images/}}

\newcommand{\nc}{\newcommand}
\nc{\RR}{\mathbb{R}}
\nc{\CC}{\mathbb{C}}
\nc{\NN}{\mathbb{N}}
\nc{\mrm}{\mathrm}
\nc{\mJ}{\mrm{J}}
\nc{\Tri}{\mathcal{T}}
\nc{\Edge}{\mathcal{E}}
\nc{\Skel}{\Gamma}
\nc{\Vh}{\mrm{V}}
\nc{\mbVh}{\mathbb{V}}
\nc{\mbVs}{\mrm{V}_{\textsc{s}}}
\nc{\bx}{\boldsymbol{x}}
\nc{\by}{\boldsymbol{y}}
\nc{\bz}{\boldsymbol{z}}
\nc{\bu}{\boldsymbol{u}}
\nc{\bn}{\boldsymbol{n}}
\nc{\bv}{\boldsymbol{v}}
\nc{\bw}{\boldsymbol{w}}
\nc{\br}{\boldsymbol{r}}

\nc{\bff}{\boldsymbol{f}}
\nc{\mbf}{\mathbf{f}}
\nc{\bg}{\boldsymbol{g}}
\nc{\bvphi}{\boldsymbol{\varphi}}
\nc{\ba}{\boldsymbol{a}}
\nc{\bb}{\boldsymbol{b}}
\nc{\bp}{\boldsymbol{p}}
\nc{\bq}{\boldsymbol{q}}
\nc{\bF}{\boldsymbol{F}}

\nc{\balpha}{\boldsymbol{\alpha}}
\nc{\bbeta} {\boldsymbol{\beta}}

\nc{\calF}{\mathcal{F}}
\nc{\calNh}{\mathcal{N}_h}

\nc{\llpar}{(\!(}
\nc{\rrpar}{)\!)}
\nc{\lbr}{\lbrack}
\nc{\rbr}{\rbrack}
\nc{\bfS}{\mathbf{S}}
\nc{\bfA}{\mathbf{A}}
\nc{\bfC}{\mathbf{C}}
\nc{\bfM}{\mathbf{M}}
\nc{\bfB}{\mathbf{B}}
\nc{\bfQ}{\mathbf{Q}}
\nc{\bfR}{\mathbf{R}}
\nc{\bfT}{\mathbf{T}}
\nc{\bfP}{\mathbf{P}}
\nc{\bfD}{\mathbf{D}}

\nc{\Nsys}{\mathrm{N}}
\nc{\bfPi}{\boldsymbol{\Pi}}
\nc{\bfId}{\mathbf{I}}

\nc{\bfH}{\mathbf{H}}
\nc{\mL}{\mrm{L}}
\nc{\bfcurl}{\mathbf{curl}}

\nc{\EdgeStar}{\Sigma}
\nc{\rg}{\operatorname{range}}

\nc{\bE}{\boldsymbol{E}}
\nc{\bH}{\boldsymbol{H}}
\nc{\bJ}{\boldsymbol{J}}

\nc{\In}{\textsc{I}}
\nc{\Bd}{\Skel}

\newtheorem{lem}{Lemma}
\newtheorem{prop}{Proposition}

\algnewcommand\algorithmicinput{\textbf{Input:}}
\algnewcommand\algorithmicoutput{\textbf{Output:}}
\algnewcommand\Input{\item[\algorithmicinput]}%
\algnewcommand\Output{\item[\algorithmicoutput]}%

\title{Nonlocal Optimized Schwarz Methods for time-harmonic electromagnetics}

\author{X. Claeys\thanks{\textsc{Ljll} (Sorbonne Universit\'e-Universit\'e de Paris-CNRS), INRIA, France}
\and F. Collino\thanks{\textsc{Poems} (CNRS-INRIA-ENSTA Paris), IP Paris, France}
\and E. Parolin\thanks{Dipartimento di Matematica, Universit\`a degli Studi di Pavia, Italia, emile.parolin@unipv.it}
}

\begin{document}

\maketitle

\begin{abstract}
  We introduce a new domain decomposition strategy for time harmonic Maxwell's
  equations that is valid in the case of automatically generated subdomain
  partitions with possible presence of cross-points. The convergence of the
  algorithm is guaranteed and we present a complete analysis of the matrix form
  of the method. The method involves transmission matrices responsible for
  imposing coupling between subdomains. We discuss the choice of such matrices,
  their construction and the impact of this choice on the convergence of the
  domain decomposition algorithm.  Numerical results and algorithms are
  provided.
\end{abstract}


\section*{Introduction}

In the context of wave propagation problems, it is known since the pioneering
work of B. Despr\'es~\cite{zbMATH00863010} that impedance type transmission
conditions shall be used between subdomains to obtain convergence of
non-overlapping domain decomposition methods (DDM).
The class of such methods is often termed Optimized Schwarz Methods (OSM).
In the simplest version of the method, the impedance operator
introduced in the transmission conditions is local.
Several alternatives for such operators were advocated,
featuring both zeroth and second order (surface) differential operators.
Without being exhaustive we mention for the acoustic
setting~\cite{Piacentini1998, Gander2002, Gander2007} and for the
electromagnetic case~\cite{Collino1997, Rodriguez2006, Dolean2009, Peng2010,
  Dolean2015b,MR2970396,MR3422381,MR3847438}.
These operators are often, but not always, constructed by mimicking absorbing
boundary conditions. For this reason, it was proposed to approximate exact absorbing conditions
by means of rational fractions of second order surface differential operators.
This was done first for the Helmholtz equation~\cite{Boubendir2012} and
then~\cite{ElBouajaji2015} for the Maxwell case.
Alternatively, non-local impedance operators were advocated in order to obtain
geometric convergence of the iterative solvers in the continuous analysis
setting~\cite{MR3989867,zbMATH01488462,zbMATH07197799}.
Such a result is out of reach with local operators for which one obtains
algebraic convergence of the DDM in the best cases, see~\cite[Chap.3]{lecouvez:tel-01444540}.

The presence of so-called cross-points i.e. points where strictly more than two
subdomains meet, has been a major and ubiquitous difficulty in the design and analysis of
efficient OSM strategies. For methods using second order surface differential operators,
cross-points are associated to corners and motivate the development of compatibility
conditions to mitigate their effects~\cite{Despres2020, Despres2021,
Modave2020, Modave2020b}.
Several other treatments inspired by available strategies developed for
elliptic problems have been proposed for nodal type
discretizations~\cite{Bendali2006,Gander2016c}.
Recently the geometric convergence result of~\cite{zbMATH01488462, zbMATH07197799}
have been extended to arbitrary geometric partitions,
including partitions with cross-points~\cite{claeys2019new,claeys2020robust}.
The new approach is based on a novel operator that communicates information
globally between subdomains and replaces the standard local exchange operator
that operates pointwise on the interface.
In addition, the method, which is derived and analysed in the acoustic setting,
is proved to be uniformly stable with respect to the discretization parameter.

We extend the work of~\cite{claeys2020robust} in three directions.
First, instead of the acoustic setting, we consider the case of electromagnetic
wave propagation problems. While no convergence result for OSM applied to Maxwell problems
in such a general context is known to us, the present analysis leads to a  convergence estimate
(see coercivity property in Proposition \ref{prop:well-posedness}) valid in the case of
heterogeneous media and general non-overlapping partitions, including
the possibility of cross-points. In the case of diagonal impedance, this yields a
new result on the pre-existing DDM strategy of Despr\'es applied to harmonic Maxwell's equations. 

Second, starting from the original undecomposed linear system,
we perform the complete derivation of the domain decomposition method
and its analysis using only matrix notations. We discard considerations related to
functional analysis and only rely on finite dimensional linear algebra and matrix
calculus so as to ease the understanding of our method in the perspective of actual
implementation. In particular Section~\ref{sec:algorithms} provides explicit algorithms.

Third, we describe a new treatment of transmission conditions that possibly
lead to extended interfaces, see Figure \ref{FigSkeleton}. In this new approach, the external
boundary of the computational domain is not necessarily part of the skeleton where transmission
conditions are imposed, which is new and computationally more optimal compared
to~\cite{claeys2019new,claeys2020robust}. 

The outline of the present contribution is as follows.
In Section~\ref{sec:geometry} we introduce several definitions and the main
notations.
In Section~\ref{sec:OrtogonalProjection} we describe the central ingredients of
our method namely the transmission matrices, the associated orthogonal projection
and the communication matrix which concentrates the main originality of the
approach.
Subsequently, the reformulation of the original problem as a skeleton problem
common to OSM is addressed in Section~\ref{sec:reformulation} followed by the
analysis of the formulation that ends with the well-posedness and convergence
results given in Proposition~\ref{prop:well-posedness}.
Next we provide two concrete choices for the transmission matrices in
Section~\ref{sec:transmission}.
The first transmission matrix stems from a simple zeroth-order
operator corresponding to the impedance operator of Despr\'es. 
The second transmission matrix stems from a more involved non-local
operator that appears to us as one of the most robust choice.
We explain in particular how to implement efficiently the latter operator
despite its underlying non-local nature.
This is followed by Section~\ref{sec:algorithms} in which we provide the
detailed algorithms in view of practical implementation of the method.
We conclude with some numerical results in Section~\ref{sec:numerics}.
In particular, we provide a first particular test case that aims at
illustrating the need for the approach that we advocate.
Besides, we investigate the influence of several parameters:
mesh refinement, wavenumber and number of subdomains.
Finally, a more involved problem featuring heterogeneous media is provided as
evidence of the robustness of the approach.

\section{Sub-domain partitioning}\label{sec:geometry}

\subsection{Mesh and vector spaces}
\begin{figure}[htp]
  \centering
  \includegraphics[width=\textwidth]{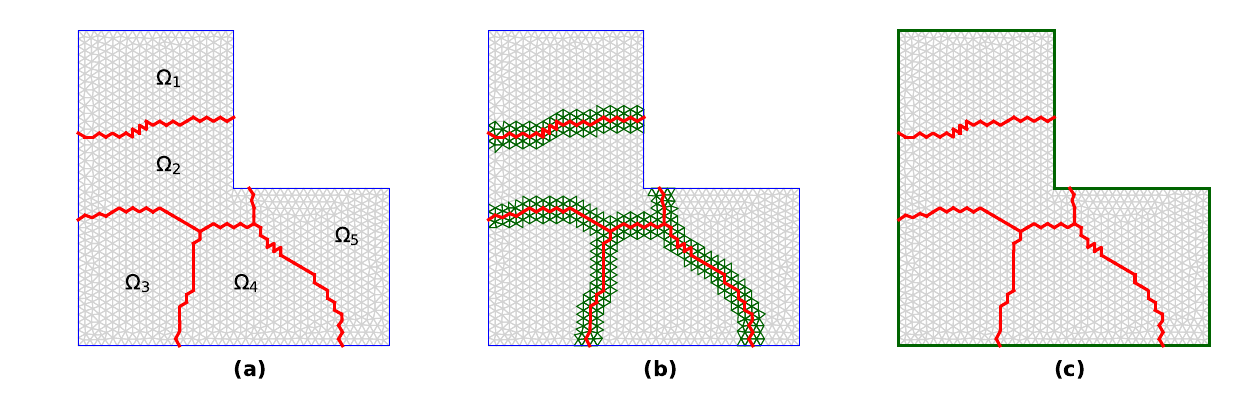}%
  \caption{A 2D sketch of extended skeleton in the case of a partition in \(5\) subdomains,
    with skeleton $\Sigma$ colored in red. Among many possibilities, the extended skeleton
    (edges colored in green and red) can be reduced to $\Sigma$ (a), it
    may consist in a thick neighborhood of $\Sigma$ (b), it may include both
    $\Sigma$ and the external boundary of the computational domain (c),
    or even a combination of the latter two sub-cases.}
  \label{FigSkeleton}
\end{figure}

We consider a (bounded) polyhedral computational domain $\Omega\subset \RR^{3}$  and a regular simplicial
triangulation $\Tri(\Omega)$ of the domain $\overline{\Omega} = \cup_{\tau\in \Tri(\Omega)}\overline{\tau}$. 
We consider a non-overlapping domain decomposition $\overline{\Omega} = \overline{\Omega}_1\cup\dots\cup
\overline{\Omega}_\mJ $ of the computational domain that is conforming with respect to the triangulation
i.e. we have the following additional properties 
\begin{equation}
  \begin{aligned}
    i)  & \quad \Omega_j\cap \Omega_k = \emptyset\;\text{if}\;j\neq k\\
    ii) & \quad \text{each $\Omega_j$ is resolved by $\Tri(\Omega)$.}
  \end{aligned}
\end{equation}
In the sequel, we shall denote $\Tri(\Omega_j) := \{\tau\in \Tri(\Omega), \tau\subset \Omega_j\}$
which implies in particular $\Tri(\Omega) = \Tri(\Omega_1)\cup \dots \cup\Tri(\Omega_\mJ)$.
This is the usual setting of non-overlapping substructuring domain decomposition methods.
Since, later on, N\'ed\'elec edge elements will be used, we introduce notations
for the edges of the mesh.
We shall denote $\Edge$ (resp. $\Edge_j$) the edges of the triangulation $\Tri(\Omega)$
(resp. $\Tri(\Omega_j)$). In particular we have $\Edge = \Edge_1\cup\dots\cup \Edge_\mJ$
which is a partition with overlap i.e. we a priori have $\Edge_j\cap \Edge_k\neq \emptyset$ if
$\Omega_j$ and $\Omega_k$ are neighboring subdomains. This leads to considering
\begin{equation}\label{DefThinSkeleton}
  \EdgeStar:= \mathop{\cup}_{\substack{1\leq j<k\leq\mJ\\ j\neq k}}\Edge_j\cap \Edge_k
\end{equation}
The edges of $\EdgeStar$ provide a triangulation of what is usually called the
\emph{skeleton} in domain decomposition literature.
Finally, we assume to have chosen a particular collection \(\Skel\) of edges
satisfying the following property 
\begin{equation}\label{DefThickSkeleton}
    \EdgeStar \subset \Skel\subset \Edge.
\end{equation}
and we set \(\Skel_j:=\Skel\cap \Edge_j\), for \(j=1,\dots,\mJ\).
We will refer to \(\Gamma\) as the \emph{extended skeleton}.
The choice of $\Skel$ satisfying the condition above may be arbitrary. 
Of course $\Skel=\EdgeStar$ (see Figure~\ref{FigSkeleton}(a))
is one possible choice among many\footnote{At first reading, one can safely
assume that $\Skel=\EdgeStar$ for simplicity.},
but the forthcoming analysis is not restricted to this sole possibility.
In practice $\Skel$ may also be chosen as a set of edges surrounding the
interfaces of the decomposition (see Figure~\ref{FigSkeleton}(b))
but we did not explore this possibility further.
An alternative is to included in \(\Skel\) the edges with multiplicity one
that belong to the physical boundary of \(\Omega\)
(see Figure~\ref{FigSkeleton}(c)).
This might have an interest if the boundary condition does not impose some more
regularity than the natural one on the associated trace.
In particular, we use this feature in some of our numerical experiments.

We also need to introduce vector spaces attached to the sets we just defined.
In the forthcoming analysis, if $\mathcal{F}$ is any finite set, we shall denote $\Vh(\mathcal{F})$
as the vector space of complex valued tuples indexed by $\mathcal{F}$ equipped with its canonical
euclidean scalar product i.e.
\begin{equation*}
    \Vh(\mathcal{F}):=\{ \bx = (x_f)_{f\in \mathcal{F}}, x_f\in \CC\}.
\end{equation*}
Elements of $\Vh(\mathcal{F})$ are tuples that may be equivalently regarded as maps
$\bx:f\mapsto x_f$ from $\mathcal{F}$ into $\CC$. Any linear map from one such space to another
$\bfM:\Vh(\calF_1)\to \Vh(\calF_2)$ is nothing but a matrix $\bfM =
(\bfM_{e,f})\in \CC^{\#\calF_2\times\#\calF_1}$ where we denoted by \(\#\calF\)
the cardinal of the set \(\calF\).
Following these notations, we can form in particular local spaces
$\Vh(\Edge_j)$ and $\Vh(\Skel_j)$ attached to each subdomain. We shall also consider
cartesian products of these spaces: for $\calF = \Edge,\Skel$ we set 
\begin{equation}
    \calF_\oplus:= \calF_1\times \dots\times \calF_\mJ\qquad\text{and}\qquad
    \Vh(\calF_\oplus) := \Vh(\calF_1)\times\dots\times \Vh(\calF_\mJ).
\end{equation}
We shall refer to $\Vh(\Skel_\oplus)$ as the \emph{multi-trace} space. This will be the space where we shall
write our final reformulation of the boundary value problem to be solved.
Let us emphasize that we use the term 'trace' even in the case where the
skeleton is extended.
Our final numerical method will take the form of a linear system posed
in $\Vh(\Skel_\oplus)$.
The size of the final matrix will then be
$\mrm{dim}\Vh(\Skel_{\oplus}) = \#\Skel_{\oplus}
= \# \Skel_1+\dots +\# \Skel_\mJ$.
We emphasize that $\#\Skel_{\oplus} > \#\Gamma$
because of overlapping between local
edge sets \(\Skel_{j}\) i.e. $\EdgeStar$
defined by \eqref{DefThinSkeleton} is a priori non-trivial.

\subsection{Restriction matrices}\label{RestrictionMatrices}
As is standard in domain decomposition, we need to introduce
restriction matrices. First we introduce $\bfR_j:\Vh(\Edge)\to \Vh(\Edge_j)$
i.e. $\bfR_j\in\CC^{\#\Edge_j\times \#\Edge}$.
These restriction matrices are collected in a global matrix (that is not a restriction
matrix)
$\bfR:\Vh(\Edge)\to \Vh(\Edge_\oplus)$
defined as follows
\begin{equation}
    \bfR^\top = \lbr \bfR_1^\top,\dots,\bfR_\mJ^\top\rbr
    \qquad\text{with}\qquad
    \bfR_j(\bx) := (x_e)_{e\in\Edge_j}\;\; \text{for}\;\;\bx = (x_e)_{e\in \Edge},
\end{equation}
where ``$\top$'' stands for the usual matrix transpose. 
The matrix $\bfR$ is a boolean matrix, by which we mean that its entries can
only take the values $0$ and $1$.
Since $\Edge = \Edge_1\cup\dots\cup \Edge_\mJ$ the matrix $\bfR$ is injective $\ker(\bfR) = \{0\}$,
but it is not surjective in general, which systematically occurs whenever
$\EdgeStar \neq \emptyset$.
Hence $\Vh(\Edge)$ is isomorphic to the range of the matrix $\bfR$ which
we shall denote by
\begin{equation}
  \mbVs(\Edge) := \rg(\bfR)\subset \Vh(\Edge_\oplus).
\end{equation}
Next we introduce similar restriction matrices
associated to the extended skeleton $\bfQ_j:\Vh(\Skel)\to \Vh(\Skel_j)$ i.e.
$\bfQ_j\in\CC^{\#\Skel_j\times \#\Skel}$.
These matrices are also collected in a global matrix (that is not a restriction
matrix) $\bfQ:\Vh(\Skel)\to \Vh(\Skel_\oplus)$ i.e.
$\bfQ\in \CC^{\#\Skel_{\oplus}\times \#\Skel}$ defined as follows
\begin{equation}\label{DefMatrixQ}
    \bfQ^\top = \lbr \bfQ_1^\top,\dots,\bfQ_\mJ^\top\rbr
    \qquad\text{with}\qquad
    \bfQ_j(\bx) := (x_e)_{e\in\Skel_j}\;\; \text{for}\;\;\bx = (x_e)_{e\in \Skel}.
\end{equation}
The matrix $\bfQ$ is also a boolean matrix and there is only one
single nonzero entry on each line.
Similarly as for \(\bfR\), the matrix \(\bfQ\) is not surjective,
but from the covering property $\Skel = \Skel_1\cup \dots\cup \Skel_\mJ$, it
follows that $\ker(\bfQ) = \{0\}$.
Hence $\Vh(\Skel)$ is isomorphic to the range of the matrix $\bfQ$ which
we shall denote by
\begin{equation}
  \mbVs(\Skel) := \rg(\bfQ)\subset \Vh(\Skel_\oplus).
\end{equation}
The space above will be referred to as the \emph{single-trace} space. It consists in those subdomain
boundary tuples that match across interfaces. Characterization of this space will be pivotal
in the forthcoming analysis. 
\begin{figure}[h]
  \centering
  \begin{subfigure}[c]{0.4\textwidth}
    \begin{tikzpicture}
      \draw[] (0, 1) node[]{\(\Vh(\Edge)\)};
      \draw[] (0.6, 1) -- (1.75, 1);
      \draw[] (2, 1) node[]{\(\bfR_{j}\)};
      \draw[->,>=latex] (2.25, 1) -- (3.4, 1);
      \draw[] (4, 1) node[]{\(\Vh(\Edge_{j})\)};
      \draw[] (4, 0.75) -- (4, 0.25);
      \draw[] (4, 0) node[]{\(\bfB_{j}\)};
      \draw[->,>=latex] (4,-0.25) -- (4,-0.75);
      \draw[] (0,-1) node[]{\(\Vh(\Skel)\)};
      \draw[] (0.6,-1) -- (1.75,-1);
      \draw[] (2,-1) node[]{\(\bfQ_{j}\)};
      \draw[->,>=latex] (2.25,-1) -- (3.4,-1);
      \draw[] (4,-1) node[]{\(\Vh(\Skel_{j})\)};
    \end{tikzpicture}
  \end{subfigure}
  \begin{subfigure}[c]{0.49\textwidth}
    \begin{tikzpicture}
      \draw[] (0, 1) node[]{\(\Vh(\Edge)\)};
      \draw[] (0.6, 1) -- (1.75, 1);
      \draw[] (2, 1) node[]{\(\bfR\)};
      \draw[->,>=latex] (2.25, 1) -- (3.4, 1);
      \draw[] (3.5, 1) node[right]{\(\mbVs(\Edge) \subsetneq \Vh(\Edge_{\oplus})\)};
      \draw[] (5.75, 0.75) -- (5.75, 0.25);
      \draw[] (5.75, 0) node[]{\(\bfB\)};
      \draw[->,>=latex] (5.75,-0.25) -- (5.75,-0.75);
      \draw[] (0,-1) node[]{\(\Vh(\Skel)\)};
      \draw[] (0.6,-1) -- (1.75,-1);
      \draw[] (2,-1) node[]{\(\bfQ\)};
      \draw[->,>=latex] (2.25,-1) -- (3.4,-1);
      \draw[] (3.5,-1) node[right]{\(\mbVs(\Skel) \subsetneq \Vh(\Skel_{\oplus})\)};
    \end{tikzpicture}
  \end{subfigure}
  \caption{Sketch of the mapping properties.
  The arrows denote surjective maps.
  }\label{fig:mappings}
\end{figure}

Next we also need to introduce trace matrices that map from the interior of
subdomains to the extended skeleton.
We introduce matrices $\bfB_j:\Vh(\Edge_j)\to\Vh(\Skel_j)$ i.e.
$\bfB_j\in \CC^{\#\Skel_j\times \#\Edge_j}$ and $\bfB:\Vh(\Edge_\oplus)\to \Vh(\Skel_\oplus)$ as follows:
\begin{equation}
    \bfB := \mrm{diag}(\bfB_1,\dots , \bfB_\mJ)
    \qquad\text{with}\qquad
    \bfB_j(\bx) := (x_e)_{e\in \Skel_j}\;\; \text{for}\;\;\bx = (x_e)_{e\in \Edge_j}.
\end{equation}
These are boolean matrices and they have only one non-zero entry per line.
The matrices $\bfB_j^{\top}$ provide a lifting from $\Vh(\Skel_j)$
into $\Vh(\Edge_j)$. In particular we have $\bfB\bfB^\top = \bfId$ and thus
$\bfB^\top\bfB$ is a projection whose action consists in cancelling those
components that are not located on $\Skel_1\times\dots\times\Skel_\mJ$.

\section{Orthogonal projection onto single traces}\label{sec:OrtogonalProjection}

\subsection{Characterizations of the single-trace space}

We start with a simple characterization of the space of single traces.

\begin{lem}\label{GlueingLemma}\quad\\
  A tuple of local subdomain contributions \(\bu\in\Vh(\Edge_\oplus)\) stems
  from a single global vector in \(\Vh(\Edge)\)
  if and only if its (interior) traces at the boundary of subdomains match at
  all interfaces.
  This is summarized as
  \begin{equation*}
    \forall \bu\in\Vh(\Edge_{\oplus}),\qquad
    \bu\in \rg(\bfR)\iff\bfB(\bu)\in \rg(\bfQ).
  \end{equation*}
\end{lem}
\noindent \textbf{Proof:}

Take an arbitrary $\bx = (\bx_1,\dots,\bx_\mJ)\in\Vh(\Edge_{\oplus})$ with $\bx_j = (x_{j,e})_{e\in \Edge_j}\in \Vh(\Edge_j)$.
Assume first that $\bx = \bfR(\by)$ for some $\by = (y_e)_{e\in\Edge}\in\Vh(\Edge)$, which writes 
$x_{j,e} = y_e$ for all $j = 1\dots \mJ$ and all $e\in \Edge_j$. Since $\Skel_j\subset \Edge_j$,
we have in particular $x_{j,e} = y_e$ for all $j = 1\dots \mJ$ and all $e\in \Skel_j$ which is
equivalent to $\bfB(\bx) = \bfQ(\bz)$ where $\bz = (y_e)_{e\in \Skel}\in \Vh(\Skel)$ i.e. $\bfB(\bx)\in
\rg(\bfQ)=\mbVs(\Skel)$.

Now assume that $\bx = (\bx_1,\dots,\bx_\mJ)\in \Vh(\Edge_\oplus)$ is such that $\bfB(\bx)\in \mbVs(\Skel)=\rg(\bfQ)$.
As a consequence there exists $\bz = (z_e)_{e\in \Skel}\in \Vh(\Skel)$ satisfying
$x_{j,e} = z_e$ for all $j = 1\dots \mJ$ and all $e\in \Skel_j = \Edge_j\cap \Skel$.
Next observe that $\Edge = \Skel\cup (\Edge_1\setminus\Skel)\cup \dots \cup  (\Edge_\mJ\setminus\Skel)$
is a disjoint union due to $\EdgeStar\subset \Skel$, see \eqref{DefThinSkeleton} and \eqref{DefThickSkeleton}.
This means that, for any $e\in \Edge$, either $e\in\Skel$, or there exists a unique $j$ such that
$e\in \Edge_j\setminus\Skel$. As a consequence we can define $\by = (y_e)_{e\in\Edge}\in \Vh(\Edge)$ by
$y_e = z_e$ if $e\in\Skel$ and $y_e = x_{j,e}$ if $e\in\Edge_j\setminus\Skel$. Because $x_{j,e} = z_e$ on
$\Edge_j\cap\Skel$, we conclude that $y_e = x_{j,e}$ for all $e\in \Edge_j$ and all $j=1\dots \mJ$,
which is equivalent to $\bx = \bfR(\by)$. \hfill$\Box$

\quad\\
The single trace space consists in those tuples of boundary traces that match at interfaces.
It yields a criterion on boundary traces for determining whenever a tuple of subdomain
contributions stems from a common global vector.

We will now discuss a more effective characterization of the space of single traces.
Instead of using pointwise constraints to ensure that a multitrace is a single trace,
we rely on a more general characterization using a projection. The idea rests on the
use of the following Lemma which is a direct consequence of Lemma \ref{GlueingLemma}.
\begin{lem}\label{GlueingLemma2}\quad\\
If $ \bfP:\Vh(\Skel_\oplus) \to\Vh(\Skel_\oplus)$ is any projection onto the single traces space
i.e. $\bfP^2=\bfP$ and $\rg(\bfP)=\rg(\bfQ)$, then
for any $\bu\in \Vh(\Edge_\oplus)$ we have 
$\bu\in \rg(\bfR)\iff (\bfId - \bfP)\bfB(\bu) = 0$.
\end{lem}

\noindent 
This observation points toward new ways to impose transmission conditions through interfaces.
This characterization of transmission conditions is the original point of our approach. 

The construction of appropriate instances of projection \(\bfP\) is not a
difficult task.
We first explain the simplest of those projections.
From \(\ker\bfQ=\{0\}\) we deduce that \(\bfQ\) admits a left pseudo-inverse
\(\bfQ^\dagger = \left(\bfQ^\top\bfQ\right)^{-1}\bfQ^\top\). This
operator can be computed explicitly since \(\bfQ^\top\bfQ\)
is diagonal.
We obtain that \(\bfQ^\dagger\bfQ = \bfId\) hence \(\bfQ^{\dagger}\) is in fact
a left inverse for \(\bfQ\).
Besides, \(\bfQ\bfQ^\dagger\) is a projection in \(\Vh(\Skel_{\oplus})\) which
is orthogonal with respect to the Euclidean scalar product and its range is
\(\mbVs(\Skel)\).
We explain the construction of other appropriate projectors (which are
orthogonal for different scalar products) in the next paragraph.

\subsection{Transmission matrices}
First, for each subdomain, we need to define the so-called 
\emph{local transmission matrices}
$\bfT_j:\Vh(\Skel_j)\to \Vh(\Skel_j)$.
Each (real-valued) $\bfT_j$ is assumed symmetric positive definite (SPD) which
is equivalent to imposing
\begin{equation}\label{AssumLocalScalProd}
    (\bx,\by)_{\bfT_j}:=\bx^\top\bfT_j\,\overline{\by}\;\;\text{for}\;\;\bx,\by\in \Vh(\Skel_j)
    \quad
    \text{is a scalar product over $\Vh(\Skel_j)$.}
\end{equation}
The norm associated with this scalar product will be denoted $\Vert \bx\Vert_{\bfT_j}^2:= (\bx,\bx)_{\bfT_j}$.
The domain decomposition strategy we are going to describe applies for any choice of
local transmission matrix $\bfT_j$ as long as they satisfy \eqref{AssumLocalScalProd}, and transmission matrices might be
regarded as parameters of the method we propose here.
In particular, this implies that $\bfT_j$ must be an invertible matrix.
How to choose properly such matrices depends on functional analysis
considerations that are discussed in \cite{claeys2020robust}, this choice
having an impact on both the speed of convergence and the computational cost of
our algorithms. 

Gathering local contributions into a single block diagonal matrix, we form a
\emph{global transmission matrix} \(\bfT\) acting on the multi-trace space
\begin{equation}\label{GlobalImpedanceScalarProd}
    \bfT = \mrm{diag}(\bfT_1,\dots,\bfT_\mJ)
    \qquad \text{and} \qquad
     (\bx,\by)_{\bfT} :=\bx^\top\bfT\,\overline{\by}
    = (\bx_1,\by_1)_{\bfT_1} + \dots +  (\bx_\mJ,\by_\mJ)_{\bfT_\mJ},
\end{equation}
for $\bx = (\bx_1,\dots,\bx_\mJ),\by =(\by_1,\dots,\by_\mJ)\in \Vh(\Skel_{\oplus})$
with $\bx_j,\by_j\in\Vh(\Skel_j)$. Consistently we
shall define the norm attached to this scalar product by $\Vert \bx\Vert_{\bfT}^2 = (\bx,\bx)_{\bfT}$. Clearly 
the block-diagonal matrix $\bfT$ induces a scalar product over $\Vh(\Skel_\oplus)$ and is thus invertible.

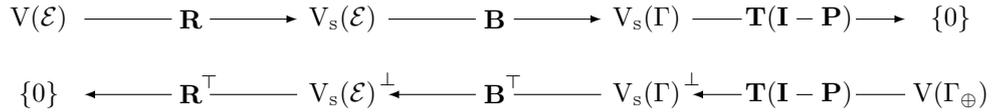
\begin{figure}[h]
  \centering
  \begin{tikzpicture}
    \draw[] (0, 0) node[]{\(\Vh(\Edge)\)};
    \draw[] (0.6, 0) -- (1.75, 0);
    \draw[] (2, 0) node[]{\(\bfR\)};
    \draw[->,>=latex] (2.25, 0) -- (3.4, 0);
    \draw[] (4, 0) node[]{\(\mbVs(\Edge)\)};
    \draw[] (4.6, 0) -- (5.75, 0);
    \draw[] (6, 0) node[]{\(\bfB\)};
    \draw[->,>=latex] (6.25, 0) -- (7.4, 0);
    \draw[] (8, 0) node[]{\(\mbVs(\Skel)\)};
    \draw[] (8.6, 0) -- (9.25, 0);
    \draw[] (10, 0) node[]{\(\bfT(\bfId-\bfP)\)};
    \draw[->,>=latex] (10.75, 0) -- (11.4, 0);
    \draw[] (12, 0) node[]{\(\{0\}\)};
    \draw[] (0, -1) node[]{\(\{0\}\)};
    \draw[<-,>=latex] (0.6, -1) -- (1.75, -1);
    \draw[] (2, -1) node[]{\(\bfR\)};
    \draw[] (2.25, -0.9) node[]{\(^{\top}\)};
    \draw[] (2.25, -1) -- (3.4, -1);
    \draw[] (4, -1) node[]{\(\mbVs(\Edge)\)};
    \draw[] (4.6, -0.85) node[]{\(^{\perp}\)};
    \draw[<-,>=latex] (4.6, -1) -- (5.75, -1);
    \draw[] (6, -1) node[]{\(\bfB\)};
    \draw[] (6.25, -0.9) node[]{\(^{\top}\)};
    \draw[] (6.25, -1) -- (7.4, -1);
    \draw[] (8, -1) node[]{\(\mbVs(\Skel)\)};
    \draw[] (8.6, -0.85) node[]{\(^{\perp}\)};
    \draw[<-,>=latex] (8.6, -1) -- (9.25, -1);
    \draw[] (10, -1) node[]{\(\bfT(\bfId-\bfP)\)};
    \draw[] (10.75, -1) -- (11.4, -1);
    \draw[] (12, -1) node[]{\(\Vh(\Skel_{\oplus})\)};
  \end{tikzpicture}
  \caption{Sketch of the mapping properties.
  The arrows denote surjective maps.
  The orthogonal complement \(^{\perp}\) is understood in the Euclidian sense.
  }\label{fig:mappings_proj}
\end{figure}

\quad\\
The forthcoming analysis will heavily rely on the projection matrix
$\bfP\in\CC^{\#\Skel_{\oplus} \times\#\Skel_{\oplus}},
  \bfP:\Vh(\Skel_\oplus)\to\mbVs(\Skel) \subset \Vh(\Skel_\oplus)$ that
is $\bfT$-orthogonal i.e. orthogonal with respect to the scalar product \eqref{GlobalImpedanceScalarProd}.
It is defined by
\begin{equation}\label{VariationnalProjection}
  \begin{aligned}
    \bfP := \bfQ(\bfQ^\top\bfT\bfQ)^{-1}\bfQ^\top \bfT
    \qquad \text{where} \qquad
    &\bfQ^\top\bfT\bfQ = \bfQ_1^\top\bfT_1\bfQ_1 + \dots + \bfQ_\mJ^\top\bfT_\mJ\bfQ_\mJ\\
    &\bfQ^\top \bfT = \lbr\bfQ_1^\top\bfT_1,\dots, \bfQ_\mJ^\top\bfT_\mJ\rbr
  \end{aligned}
\end{equation}
It appears obvious from the definition above that
\(\left[\bfT(\bfId-\bfP)\right]^{\top} = \bfT(\bfId-\bfP)\).
We state now the counterpart of Lemma~\ref{GlueingLemma2} for a
\(\bfT\)-orthogonal projection \(\bfP\).
The identities appearing in Lemma~\ref{GlueingLemma3} are represented in
Figure~\ref{fig:mappings_proj}.

\begin{lem}\label{GlueingLemma3}\quad\\
  If $ \bfP:\Vh(\Skel_\oplus) \to\Vh(\Skel_\oplus)$ is a \(\bfT\)-orthogonal
  projection onto the single traces space i.e. $\bfP^2=\bfP$,
  $\bfP$ is self-adjoint with respect to the scalar product induced by \(\bfT\)
  and $\rg(\bfP) = \mbVs(\Skel)$, then
  \begin{equation}\label{LiftingLemma}
      \rg(\bfR) = \ker(\bfT(\bfId-\bfP)\bfB)
      \qquad\text{and}\qquad
      \ker(\bfR^{\top}) = \rg(\bfB^{\top}\bfT(\bfId-\bfP)).
  \end{equation}
\end{lem}

\subsection{Projecting a multiple trace in practice}\label{sec:proj_in_practice}
Formula \eqref{VariationnalProjection} involves the inverse matrix $(\bfQ^\top\bfT\bfQ)^{-1}$ which
indicates that computing the action of the orthogonal projection $\bfP:\Vh(\Gamma_\oplus)\to\Vh(\Gamma_\oplus)$
requires the solution to an auxiliary linear system associated to the matrix $\bfQ^\top\bfT\bfQ:
\Vh(\Gamma)\to\Vh(\Gamma)$.
In practice, following the above formula for the projection matrix, the image
$\bfP(\bx)$ of any element $\bx\in \Vh(\Skel_\oplus)$
can be computed as follows
\begin{equation}\label{FormulaProjection2}
  \begin{aligned}
    & \bw = \bfP(\bx)\;\iff\;
    \bw = \bfQ(\by)\\
    &\text{where}\;\;
        \by\in \Vh(\Gamma)\;\;\text{solves}\;\;
        (\bfQ^\top\bfT\bfQ)(\by) = \bfQ^\top\bfT\bx.
  \end{aligned}
\end{equation}
This linear system is not a priori block-diagonal and, in general, will not be. It should be interpreted
as a non-local operator. Despite the nonlocality of $\bfP$ and the fact that \eqref{VariationnalProjection} only provides
an implicit definition of \(\bfP\), what really matters is fast evaluation of $\bx\mapsto \bfP(\bx)$. The bottleneck here is
of course the solution to the $\#\Gamma\times\#\Gamma$ linear system \eqref{FormulaProjection2}. This requires an
efficient solution strategy for computing \eqref{VariationnalProjection}, which involves an SPD problem.
Current literature
already provides many powerful techniques for solving such problems including
adaptive multigrids, see e.g.~\cite{zbMATH03874483}, or two-level substructuring domain
decomposition method \cite[Chap.4-6]{zbMATH02113718}.
A possible solution strategy for treating this linear system relies on the
Neumann-Neumann algorithm. Observe from Definition \eqref{DefMatrixQ} that the matrix
$\bfD := (\bfQ^\top\bfQ)^{-1}:\Vh(\Gamma)\to \Vh(\Gamma)$ is diagonal $\bfD = \mrm{diag}_{e\in \Gamma}(1/d_e)$
where $d_e =\#\{j \in \{1,\dots, \mJ\}, \;e\in \Edge_j\}$ is the number of
subdomains an $e$ belongs to. The Neumann-Neumann algorithm~\cite{zbMATH06029154,zbMATH02113718} then consists in a
preconditioned conjugate gradient solver (PCG) taking $\bfD\,\bfQ^\top \bfT^{-1}\bfQ\,\bfD$
as preconditioner. 
Let \(\bfM:= \bfD\,\bfQ^\top \bfT^{-1}\bfQ\,\bfD\), the preconditioned problem then writes 
\begin{equation}\label{NeumannNeumannStrategy}
    \bw = \bfP(\bx)\;\;\iff\;\;
    \bw = \bfQ(\by)\;\;\text{where}
    \;\;\by\in \Vh(\Gamma)\;\;\text{solves}\;\; 
    \bfM(\bfQ^\top\bfT\bfQ) (\by) = \bfM\bfQ^\top \bfT(\bx).
\end{equation}

\subsection{Communication matrix}
The projection $\bfP$ leads to the definition of a so-called \textit{communication matrix}
$\bfPi\in \CC^{\#\Skel_{\oplus}\times\#\Skel_{\oplus}}$ defined as the matrix of the orthogonal symmetry with respect
to $\mbVs(\Skel)$ i.e.
\begin{equation}\label{ExchangeMatrixDef}
    \bfPi:=2\bfP-\bfId
    \;\;\text{so that}\;\;
    \bfP = (\bfId+\bfPi)/2.
\end{equation}
Observe that $\bfPi = \bfP - (\bfId - \bfP)$ and that $\bfId - \bfP = (\bfId-\bfPi)/2$ which is the $\bfT$-orthogonal
projection with $\mbVs(\Skel)$ as kernel. The communication matrix satisfies a few elementary yet important
properties that are summarized in the next lemma.
\begin{lem}\label{PropertiesExchangeMatrix}\quad\\
  The communication matrix $\bfPi$ defined by \eqref{ExchangeMatrixDef} is a $\bfT$-isometric
  involution i.e. $\bfPi^2 = \bfId$ and $\Vert \bfPi(\bx)\Vert_\bfT = \Vert \bx\Vert_\bfT$ for all
  $\bx\in \Vh(\Skel_\oplus)$.
\end{lem}
\noindent \textbf{Proof:}

We have $\bfP^2 = \bfP$ since $\bfP$ is a projection by construction, so that $\bfPi^2 = 4\bfP^2 - 4\bfP + \bfId = \bfId$.
On the other hand, from  \eqref{VariationnalProjection} we conclude that
$\bfP^{\top}\bfT\bfP = \bfT\bfP$ hence  $\Vert \bfP(\bx)\Vert_\bfT^2 =
(\bfP(\bx),\bfP(\bx))_\bfT = (\bx,\bfP(\bx))_\bfT = \Re e\{(\bx,\bfP(\bx))_\bfT\}$. As a
consequence
\begin{equation*}
  \begin{aligned}
    \Vert \bfPi(\bx)\Vert_\bfT^2
    & = \textcolor{white}{4}\Vert 2\bfP(\bx) - \bx\Vert_\bfT^2 = (2\bfP(\bx) - \bx,2\bfP(\bx) - \bx)_\bfT\\
    & = 4\Vert \bfP(\bx)\Vert_\bfT^2 - 4\Re e\{ (\bx,\bfP(\bx))_\bfT\} + \Vert \bx\Vert_\bfT^2 = \Vert \bx\Vert_\bfT^2.
  \end{aligned}
\end{equation*}
\hfill $\Box$

\subsection{Explicit expressions}\label{ExplicitExpression}
Although the projection and communication matrices $\bfP$ and $\bfPi$ are non-local in general,
there are cases where they get localized. There are choices of $\bfT$ for which it is
possible to exhibit an explicit expression
for the matrix $(\bfQ^\top\bfT\bfQ)^{-1}\bfQ^\top\bfT(\bu)$. 
A first simple example is the case of a scalar transmission matrix, namely
\(\bfT=a\mathbf{I}\) (with \(a>0\)) for which we immediately get 
$\bfP = \bfQ(\bfQ^\top\bfQ)^{-1}\bfQ^\top$. 
Remarkably, the projection is independent of \(a\) (hence of \(\bfT\)).

Next, we give another example that is a generalization of the previous simple
case to some diagonal matrices.
This is a fundamental particular case since it corresponds to the 
overwhelming majority of domain decomposition methods where the exchange of
information between adjacent subdomains simply consists in swapping data through
their common interface.
Deviations from this case mainly include special treatments for geometries with
cross-points.
For each $e\in \Gamma$ set $\Upsilon(e):=\{j\in \{1,\dots,\mJ\}, e\in \Gamma_j\}$
and \(d_{e} = \#\Upsilon(e)\).
Then, for any subset $\Upsilon\in\mathscr{P}(\{1,\dots,\mJ\})$ where $\mathscr{P}(E)$ refers
to the subsets of $E$, denote $\Gamma_\Upsilon := \{e\in \Gamma, \Upsilon(e) = \Upsilon\}$.
The collection of $\Gamma_\Upsilon$ yields a disjoint partition of $\Gamma$ associated to
the equivalence relation $e\sim e'\iff \Upsilon(e) = \Upsilon(e')$, see
\cite[\S 2.5.1]{zbMATH06029154}. In particular
$\Gamma_\Upsilon\cap \Gamma_{\Upsilon'} = \emptyset$ if $\Upsilon \neq \Upsilon'$. Next
consider scalar products defined through the symmetric positive definite matrices
$\bfT_\Upsilon: \Vh(\Gamma_\Upsilon)\to \Vh(\Gamma_\Upsilon)$, and assume that 
each local transmission matrix $\bfT_j:\Vh(\Gamma_j)\to \Vh(\Gamma_j)$ satisfies
\begin{equation}\label{BlockDiagonalForm}
  \begin{aligned}
    & (\bfT_j)_{e,e'} = 0 &&\text{if}\quad\Upsilon(e)\neq \Upsilon(e')\\
    & (\bfT_j)_{e,e'} = (\bfT_\Upsilon)_{e,e'} &&\text{if}\quad\Upsilon(e) = \Upsilon(e') = \Upsilon.
  \end{aligned}
\end{equation}
This means that each local transmission matrix $\bfT_j$ is assumed block diagonal,
each block $\bfT_\Upsilon$ corresponding to one of the equivalence classes intersecting
$\Gamma_j$. With such a choice of transmission matrix, then $\bv = \bfP(\bu)$ is given
by the explicit formula
\begin{equation}\label{OrthogonalProjection}
  \begin{aligned}
    \bv = \bfP(\bu) \iff\;\;
    & \bv_{j,e} = \frac{1}{d_{e}}\sum_{k\in \Upsilon(e)} \bu_{k,e}\;\;\forall e\in \Skel\\
    & \text{where}\quad \bu = (\bu_1,\dots,\bu_\mJ)\in\Vh(\Skel_\oplus),\;\;
    \bu_j = (\bu_{j,e})_{e\in \Gamma_j},\\
    & \textcolor{white}{where}\quad \bv = (\bv_1,\dots,\bv_\mJ)\in\Vh(\Skel_\oplus),\;\;
    \bv_j = (\bu_{j,e})_{e\in \Gamma_j}.
  \end{aligned}
\end{equation}
Notice that \(\bfP\) is independent of the transmission matrices $\bfT_j$ (hence also \(\bfPi\)).
Let us examine the particular case where the domain decomposition does not
involve any cross-point and $\Skel = \EdgeStar$. Such decompositions are
sometimes referred to as "onion skin" like. Hypothesis \eqref{BlockDiagonalForm}
then means that $\bfT_j$ couples edges belonging to the same interface.
In this special case we have $d_{e} = 2$, $\forall e\in\EdgeStar$ i.e.
$\Upsilon(e) = \{j_{-}(e),j_{+}(e)\}$ where $j_{-}(e)<j_{+}(e)$ so, with the same
notation as in \eqref{OrthogonalProjection}, the orthogonal projection and the
communication matrix are fully local matrices and
are given explicitly by the formula
\begin{equation}\label{ExpressionExplicitSwapping}
  \begin{alignedat}{3}
    \bv = \bfP(\bu)  \iff \;\; &
    \bv_{j_{-}(e),e} = \bv_{j_{+}(e),e} = (\bu_{j_{-}(e),e}+ \bu_{j_{+}(e),e})/2,
    &&\qquad \forall e\in \EdgeStar,\\
    \bv = \bfPi(\bu) \iff \;\; &
    \bv_{j_{-}(e),e} = \bu_{j_{+}(e),e} \;\;\;\text{and}\;\;\;
    \bv_{j_{+}(e),e} = \bu_{j_{-}(e),e},
    &&\qquad \forall e\in \EdgeStar.
  \end{alignedat}
\end{equation}
We recover the familiar swapping of data at each interface and our
approach based on orthogonal projections is then proved to be a proper
generalization of the standard technique in domain decomposition methods.

\section{The scattering problem and its reformulation}\label{sec:reformulation}
The present contribution is concerned with the efficient solution to electromagnetic
scattering problems. Although the principles that we are going to develop apply to a
wider range of problems, for the sake of clarity, we choose a specific model problem
for explaining our method and we describe this model problem here.

\subsection{Variational problem and Galerkin approximation}
First we need to formulate a few reasonable assumptions regarding
the coefficients modelling the propagation medium. We shall assume a strictly
constant positive (angular) frequency $\omega>0$ as well as
three measurable essentially bounded functions:
the electric permittivity and the magnetic permeability
$\epsilon,\mu:\Omega\to \CC$ 
and the impedance $\eta:\partial\Omega\to \CC$.
We assume that these functions are also uniformly bounded below
i.e., there exist constants $\epsilon_\star,\mu_\star,\eta_\star>0$ such that 
$\Re e\{\epsilon(\bx)\}>\epsilon_\star, \Re e\{\mu(\bx)\}>\mu_\star$ for all
$\bx\in \Omega$ and $\Re e\{\eta(\bx)\}>\eta_\star$ for all $\bx\in
\partial\Omega$.
We also assume 
\begin{equation}\label{Dissipativity1}
    \Im m\{\epsilon(\bx)\}\geq 0, \;\; \Im m\{\mu(\bx)\}\geq 0\quad \forall \bx\in \Omega,
    \qquad\text{and}\qquad
    \Im m\{\eta(\bx)\}\geq 0\quad \forall\bx\in\partial\Omega.
\end{equation}
In the following $\bn:\partial\Omega\to \RR^3$ shall refer to the outward pointing
unit normal vector to the boundary of the computational domain.
Given a volume source term $\bJ\in L^2(\Omega^3)$ 
and a surface current $\bJ_\sigma$ 
i.e.\ a tangential vector field in $L^2(\partial \Omega)^3$ with
$\bJ_\sigma \cdot \bn =0$,
we consider the model problem: find electric and magnetic
fields $\bE,\bH\in \mL^{2}(\Omega)^3$ satisfying
\begin{equation}\label{eq:maxwell_strong_form}
  \begin{cases}
    \bfcurl(\bE) - \imath \omega\mu \bH = 0         & \text{in}\;\Omega,\\
    \bfcurl(\bH) + \imath \omega\epsilon \bE = \bJ  & \text{in}\;\Omega,\\
    \bn\times [\bE\times \bn] - \eta\,\bH\times\bn = \eta \bJ_\sigma & \text{on}\;\partial\Omega.
  \end{cases}
\end{equation}
Here of course $\bn\times [\bE\times \bn]$ is the tangential component
of the electric field on the boundary $\partial\Omega$. Eliminating the magnetic field
$\bH$, this problem can be equivalently put in variational form with the electric
field \(\bE\) as sole unknown: find $\bE\in \mathcal{W}(\Omega):=\{ \bu\in\mL^{2}(\Omega)^3, \bfcurl(\bu)\in
\mL^2(\Omega)^3, \bu\times\bn\in \mL^2(\partial\Omega)^3\}$ such that $a_\Omega(\bE,\bE') = \ell_\Omega(\bE')$ for all $\bE'\in
\mathcal{W}(\Omega)$ where
\begin{equation}
  \begin{aligned}
    a_{\Omega}(\bu,\bv):=
    & \int_{\Omega}\mu_r^{-1}\bfcurl(\bu)\cdot\bfcurl(\overline{\bv})
    - \kappa^2 \epsilon_r \bu\cdot\overline{\bv}\,d\bx
    -\imath\kappa\int_{\partial\Omega} \eta_r^{-1}(\bu\times\bn)\cdot(\overline{\bv}\times\bn) d\sigma,\\
    \ell_{\Omega}(\bv) :=
    &  \imath\kappa \int_{\Omega} \sqrt{{\mu_0^{}}{\epsilon_0^{-1}}} \bJ\cdot \overline{\bv}d\bx
    - \imath\kappa \int_{\partial\Omega} \sqrt{{\mu_0^{}}{\epsilon_0^{-1}}} \bJ_\sigma\cdot \overline{\bv}d\sigma.
  \end{aligned}
\end{equation}
Here we have introduced dimensionless and possibly varying relative
parameters (indexed by $r$),
using the constant values in the vacuum (indexed by $0$),
namely $\epsilon=\epsilon_0 \epsilon_r$, $\mu=\mu_0 \mu_r$,
$\eta=\sqrt{{\mu_0}/{\epsilon_0}}\, \eta_r$.
Besides, we denote by $\kappa=\omega \sqrt{{\mu_0}{\epsilon_0}}$
the constant wave number in the vacuum.

We consider a Galerkin discretization of this problem by means of N\'ed\'elec edge's finite
elements: find $\bE_h\in \calNh(\Omega)$ such that $a_\Omega(\bE_h,\bE_h') = \ell_\Omega(\bE_h')$
for all $\bE_h'\in \calNh(\Omega)$, with discrete variational space defined by
$\calNh(\Omega) := \{\bu\in \mathcal{W}(\Omega),
\bu\vert_{\tau}\in \mathcal{N}(\tau)\}$ and $\mathcal{N}(\tau):=\{\bvphi\vert_{\tau},
\bvphi(\bx) = \balpha+\bx\times\bbeta,\balpha,\bbeta\in\CC^3\}$. After fixing
a collection $\{\boldsymbol{t}_e\}_{e\in\Edge}$ where $\boldsymbol{t}_e\in \RR^3$
is a unit tangent vector to the edge $e$, the discrete variational 
space is decomposed according to shape functions $\calNh(\Omega)  =
\mrm{span}\{\bvphi_e(\bx),e\in\Edge\}$ where $\bvphi_e$ is the only element
$\calNh(\Omega)$ satisfying
$\int_{e}\bvphi_e(\bx)\cdot\boldsymbol{t}_e\,d\sigma(\bx) = 1$ and
$\int_{f}\bvphi_e(\bx)\cdot\boldsymbol{t}_f\,d\sigma(\bx) = 0$ for $f\in \Edge$
and $f\neq e$. We finally obtain the matrix form of the problem: noting
$\bu_\Omega = (\int_{e}\bE_h\cdot\boldsymbol{t}_e\,d\sigma)_{e\in\Edge}$, we look for
\begin{equation}\label{Formulation0}
  \begin{aligned}
    \bu_\Omega\in \Vh(\Edge)\;\;\text{such that}\;\;\bfA_{\Omega} \bu_\Omega = \bff_\Omega
    \;\;\text{where}\;\;
    \bfA_{\Omega} = (a_{\Omega}(\bvphi_f,\bvphi_e))_{e,f\in \Edge},\ 
    \bff_{\Omega} = (\ell_{\Omega}(\bvphi_e))_{e\in \Edge}.
  \end{aligned}
\end{equation}
Provided that the mesh is sufficiently fine, which we shall systematically assume thereafter,
it is a consequence of classical analysis of Maxwell's equations \cite{zbMATH01266750,zbMATH05080480,zbMATH01868846}
that Assumption \eqref{Dissipativity1} implies the well posedness of the Galerkin variational
formulation \eqref{Formulation0} and thus the invertibility of the matrix $\bfA_{\Omega}$.

\subsection{Reformulation based on domain decomposition}
Domain decomposition leads to considering restricted sesquilinear forms on each
local subdomain. Denote $\calNh(\Omega_j):=\{\bv\vert_{\Omega_j},\;\bv\in \calNh(\Omega)\}$,
and 
\begin{equation*}
  \begin{aligned}
    a_{\Omega_j}(\bu,\bv):=
    & \int_{\Omega_j}\mu_r^{-1}\bfcurl(\bu)\cdot\bfcurl(\overline{\bv})
    - \kappa^2 \epsilon_r \bu\cdot\overline{\bv}\,d\bx
    -\imath\kappa\int_{\partial\Omega_j\cap\partial\Omega} \eta_r^{-1}(\bu\times\bn)\cdot(\overline{\bv}\times\bn) d\sigma,\\
    \ell_{\Omega_j}(\bv) :=
    & \imath\kappa \int_{\Omega_j} \sqrt{{\mu_0^{}}{\epsilon_0^{-1}}} \bJ\cdot \overline{\bv}d\bx
    - \imath\kappa \int_{\partial\Omega_j\cap\partial\Omega} \sqrt{{\mu_0^{}}{\epsilon_0^{-1}}} \bJ_\sigma\cdot \overline{\bv}d\sigma.
  \end{aligned}
\end{equation*}
The corresponding local matrices and right hand sides take the 
expression
\begin{equation}\label{LocalContrib}
    \bfA_{j}   := (a_{\Omega_j}(\bvphi_f,\bvphi_e))_{e,f\in\Edge_j}
    \qquad\text{and}\qquad
    \bff_j \;\,:= (\ell_{\Omega_j}(\bvphi_e))_{e\in\Edge_j}.
\end{equation}
In particular $\bfA_{j}\in \CC^{\#\Edge_j\times\#\Edge_j}$ and $\bff_j\in \CC^{\#\Edge_j}$.
The local contributions \eqref{LocalContrib} are glued together, enforcing
continuity across interfaces, by means of the restrictions matrices $\bfR_j$
introduced in \S\ref{RestrictionMatrices}. The global linear system is then
decomposed in the following manner
\begin{equation}\label{GlobalMatrix}
  \begin{aligned}
    & \bfA_\Omega = \bfR^\top\bfA\bfR = \bfR_1^\top\bfA_1\bfR_1+\dots+\bfR_\mJ^\top\bfA_\mJ\bfR_\mJ\\
    & \bff_\Omega = \bfR^\top\bff = \bfR_1^\top\bff_1+\dots+\bfR_\mJ^\top\bff_\mJ\\
    & \text{where}\quad \bfA := \mrm{diag}(\bfA_1,\dots,\bfA_\mJ)
    \quad \text{and}\quad \bff^\top :=\lbr \bff^\top_1,\dots,\bff_\mJ^\top\rbr
  \end{aligned}
\end{equation}
so that $\bfA\in \CC^{\#\Edge_{\oplus}\times\#\Edge_{\oplus}}$
and $\bff\in \CC^{\#\Edge_{\oplus}}$.

Below we introduce several problems that are equivalent to the original
Problem~\eqref{Formulation0}. By equivalent we mean that having a solution to
one of the two problems yields a solution to the other one.

\paragraph{Reformulation 1.}
Let \(\bu_{\Omega}\) be solution to the discrete problem~\eqref{Formulation0},
namely \(\bfR^{\top}\bfA\bfR\bu_{\Omega} = \bfR^{\top}\bff\).
Introducing \(\bu = \bfR\bu_\Omega\in \Vh(\Edge_\oplus)\)
and \(\bv = \bfA\bu-\bff \in \Vh(\Edge_\oplus)\) it is immediate to see that
they are solutions to
\begin{equation}\label{Formulation1}
  \text{find}\;(\bu,\, \bv)\in \rg\bfR\times\ker\bfR^{\top}
  \;\text{such that} \qquad
  \bfA\bu - \bv = \bff.
\end{equation}
Reciprocally, if \((\bu,\, \bv)\) are solutions to~\eqref{Formulation1}
then there exists \(\bu_{\Omega}\in\Vh(\Edge)\) such that \(\bfR \bu_{\Omega} = \bu\).
Multiplying both sides by \(\bfR^{\top}\bfA\) yields
\(\bfR^{\top}\bfA\bfR \bu_{\Omega} = \bfR^{\top}\bfA\bu = \bfR^{\top}(\bv+\bff) = \bfR^{\top}\bff\)
and \(\bu_{\Omega}\) is solution to~\eqref{Formulation0}.
In fact, the global solution \(\bu_{\Omega}\) might be recovered from the
solution \(\bu\) of the problem above by the identity
$\bu_\Omega = (\bfR^\top\bfR)^{-1}\bfR^\top\bu$ which does not raise any
computational difficulty since the matrix $\bfR^\top\bfR\in
\CC^{\#\Edge\times\#\Edge}$ is diagonal i.e. 
$\bfR^\top\bfR = \mrm{diag}_{e\in \Edge}(d_e)$ where
$d_e =\#\{j \in \{1,\dots, \mJ\} \;e\in \Edge_j\}$.

\paragraph{Reformulation 2.}
 Problem~\eqref{Formulation1} is block diagonal which allows
exploiting the geometric partitioning of the domain. However the solution
spaces \(\rg(\bfR)\) and \(\ker(\bfR^{\top})\) are not convenient, and
\(\bfA\) is not always invertible. This motives changes
of unkowns. First, according to Lemma \ref{GlueingLemma3},
$\bv\in \ker\bfR^{\top}$ if and only if 
$\bv = \bfB^\top\bfT(\bp+\imath \bfB\bu)$ for some $\bp\in \Vh(\Skel_\oplus)$
satisfying $\bfP(\bp+\imath\bfB\bu)=0$. Such a $\bp$ is then unique and
given by $\bp = \bfT^{-1}\bfB\bv -\imath \bfB\bu$.
In addition, Lemma \ref{GlueingLemma3} also implies that
$\bu\in \rg\bfR$ if and only if $\bu\in \Vh(\Edge_\oplus)$
and $\bfB\bu\in \ker(\bfId-\bfP)$. Hence $(\bu,\bv)$ solves
\eqref{Formulation1} if and only if $(\bu,\bp)$
solves
\begin{equation}\label{Formulation3}
  \begin{aligned}
    \text{find}\;(\bu,\bp)\in \Vh(\Edge_\oplus)\times \Vh(\Skel_\oplus)
    \;\text{such that} \qquad
    & (\bfA-\imath\bfB^\top\bfT\bfB)\bu-\bfB^\top\bfT\bp = \bff,\\
    & (\bfId-\bfP)\bfB\bu = 0,\\
    & \bfP(\bp+\imath\bfB\bu)=0.
  \end{aligned}
\end{equation}
The next lemma establishes that the linear system in the second line above is
systematically well posed. The proof relies on Assumption \eqref{Dissipativity1}
which implies in particular that $\Im m\{a_{\Omega_j}(\bv,\bv)\}\leq 0$ for all
$\bv\in \calNh(\Omega_j)$. We deduce that
$\Im m\{\overline{\bv}^\top\bfA_{j}\bv\}\leq 0$ for all $ \bv\in \Vh(\Edge_j)$.
The sign property on the imaginary part of each local contribution \(\bfA_{j}\) naturally
transfers to the global matrix \(\bfA\). This leads to introducing a quadratic
functional $\mathcal{P}: \Vh(\Edge_\oplus)\to \RR$ associated to energy dissipation
\begin{equation}\label{DefinitionDissipatedEnergy}
    \mathcal{P}(\bv):=-\Im m\{\overline{\bv}^\top\bfA\bv\},
    \qquad
    \mathcal{P}(\bv) \geq 0,\quad\forall \bv\in \Vh(\Edge_\oplus).
\end{equation}

\begin{lem}\quad\\
  The matrix $\bfA-\imath \bfB^\top\bfT\bfB$ is invertible.
\end{lem}
\noindent \textbf{Proof:}

It suffices to show that the kernel is trivial. Pick $\bu\in \Vh(\Edge_\oplus)$ satisfying
$(\bfA-\imath \bfB^\top\bfT\bfB)\bu = 0$. This implies in particular 
$\Vert \bfB\bu\Vert_{\bfT}^2 - \Im m\{\overline{\bu}^\top\bfA\bu\} = 0$ and,
taking account of \eqref{DefinitionDissipatedEnergy}, we conclude that $\bfB\bu = 0$,
and thus $\bfA\bu = 0$. Next according to \eqref{LiftingLemma}, there exists
$\bu_\Omega\in \Vh(\Edge)$ such that $\bu = \bfR(\bu_\Omega)$ so that $\bfA\bfR(\bu_\Omega)
= 0$ and thus $\bfA_\Omega(\bu_\Omega) = \bfR^\top\bfA\bfR\bu_\Omega = 0$. Since
$\bfA_\Omega$ is invertible due to well posedness of the original wave scattering problem
\eqref{Formulation0}, this implies $\bu_\Omega = 0$ and thus $\bu = \bfR(\bu_\Omega) = 0$.
\hfill $\Box$

\quad\\
It is important to realize that, in the first equation above, the matrix $\bfA-\imath \bfB^\top\bfT\bfB
= \mrm{diag}_{j=1\dots \mJ}(\bfA_j-\imath \bfB_j^\top\bfT_j\bfB_j)\in \CC^{\#\Edge_{\oplus}\times\#\Edge_{\oplus}}$ is block diagonal
since $\bfA,\bfB$ and $\bfT$ are themselves block diagonal. 

\paragraph{Reformulation 3.}
Further elaborating on Formulation \eqref{Formulation3}, since $\bfP$ is a $\bfT$-orthogonal
projection the two equations $(\bfId-\bfP)\bfB\bu = 0$ and $\bfP(\bp+\imath\bfB\bu)=0$
are equivalent to the single statement $2\bfP(\bp+\imath\bfB\bu) + 2\imath (\bfId-\bfP)\bfB\bu = 0$.
Then taking into account that $\bfP = (\bfPi+\bfId)/2$ and $\bfId-\bfP = (\bfId-\bfPi)/2$,
the latter equation can be rewritten $\bp + \bfPi(\bp+2\imath \bfB\bu) = 0$. This shows that
\eqref{Formulation3} can be written equivalently 
\begin{equation}\label{Formulation4}
  \begin{aligned}
    \text{find}\;(\bu,\bp)\in \Vh(\Edge_\oplus)\times \Vh(\Skel_\oplus)
    \;\text{such that}\qquad
    & (\bfA-\imath \bfB^\top\bfT\bfB)\bu - \bfB^\top\bfT\bp = \bff,\\
    & \bp + \bfPi(\bp+2\imath \bfB\bu)=0.
  \end{aligned}
\end{equation}

\paragraph{Reformulation 4.}
Invertibility of the matrix $\bfA-\imath \bfB^\top\bfT\bfB$ allows to eliminate the volume unknown
$\bu\in\Vh(\Edge_\oplus)$ from Problem \eqref{Formulation4} and to reduce this problem to an equation
posed on $\Skel$ only. To achieve this, let us set
\begin{equation}\label{DefScatteringMatrix}
  \begin{aligned}
    & \bfS := \bfId + 2\imath\bfB(\bfA-\imath \bfB^\top\bfT\bfB)^{-1}\bfB^\top\bfT,\\
    & \bb:=-2\imath \bfPi\bfB(\bfA-\imath \bfB^\top\bfT\bfB)^{-1}\bff.
  \end{aligned}
\end{equation}
The matrix
$\bfS\in \CC^{\#\Edge_{\oplus}\times\#\Edge_{\oplus}}$
is commonly called \emph{scattering matrix}. It is a block diagonal
matrix
\(\bfS := \operatorname{diag}\bfS_{j}\)
with \(\bfS_{j}:=\bfId + 2\imath\bfB_{j}(\bfA_{j}-\imath \bfB_{j}^\top\bfT_{j}\bfB_{j})^{-1}\bfB_{j}^\top\bfT_{j}\)
whose inversion can thus be made fully parallel. Its definition guarantees that $\bp+2\imath
\bfB\bu = \bfS(\bp) + 2\imath \bfB(\bfA-\imath \bfB^\top\bfT\bfB)^{-1}\bff$. Plugging this into the second equation of \eqref{Formulation3}, we finally
arrive at what we shall call ``skeleton formulation'', namely 
\begin{equation}\label{SkeletonFormulation}
    \text{find}\;\bp\in \Vh(\Skel_\oplus)
    \;\text{such that}\qquad
    (\bfId + \bfPi\bfS)\bp = \bb.
\end{equation}
As mentioned above, once Equation \eqref{SkeletonFormulation} is solved,
the global volume solution can be recovered by computing
$\bu = (\bfA-\imath \bfB^\top\bfT\bfB)^{-1}(\bfB^\top\bfT\bp+\bff)$ which can
be achieved in parallel since the matrix
$\bfA-\imath \bfB^\top\bfT\bfB := \operatorname{diag} (\bfA_{j}-\imath \bfB_{j}^\top\bfT_{j}\bfB_{j})$
is subdomain-wise block diagonal.

To sum up, we have given different equivalent formulations (i.e. \eqref{Formulation1}, \eqref{Formulation3}
then \eqref{Formulation4}) of the initial problem \eqref{Formulation0} and finally obtain the skeleton formulation
\eqref{SkeletonFormulation} which is the one we propose to solve by an appropriate linear solver. This form is not new;
the equation $(\bfId + \bfPi\bfS)\bp = \bb$ with $\bfPi$ the exchange matrix defined in
\S\ref{ExplicitExpression}, Formula \eqref{ExpressionExplicitSwapping}, appears in \cite{zbMATH07197799,zbMATH01488462}
where DDM algorithm is applied to Helmholtz equation with an onion skin domain decomposition i.e. no
cross-point\footnote{With a different sign convention though, that results in
considering $-\bfS$ instead of $+\bfS$.}. Although
these previous works can easily be extended to Maxwell's equations, they can only handle interfaces with edges of
multiplicity two. Here we obtain a generalization that yields a treatment of cross-points with edges of greater
multiplicity. The price to pay is a more elaborate definition of $\bfPi$.

\subsection{Analysis of the skeleton formulation}\label{ScatteringMatrix}
Besides the communication matrix $\bfPi$, the scattering matrix $\bfS$ is
a cornerstone of Equation \eqref{SkeletonFormulation}. It models the wave
propagation phenomena within each local subdomain. We dedicate the
present section to deriving a few key properties of this matrix.
Taking account of the identity $\bfId = \bfB\bfB^\top = \bfB(\bfA-
\imath\bfB^\top\bfT\bfB)^{-1} (\bfA-\imath\bfB^\top\bfT\bfB)\bfB^\top$,
a basic simple re-arrangement in the definition \eqref{DefScatteringMatrix}
of the scattering matrix $\bfS$ yields the expression
\begin{equation}\label{FirstExpressionScattering}
\bfS = \bfB(\bfA-\imath\bfB^\top\bfT\bfB)^{-1}
(\bfA+\imath\bfB^\top\bfT\bfB)\bfB^\top.
\end{equation}
This expression can be further condensed by means of the Schur complement of the
matrix $\bfA = \mrm{diag}(\bfA_1,\dots,\bfA_\mJ)$ following the standard approach
in substructuring methods \cite[chap.2]{zbMATH06029154}, \cite[chap.4]{zbMATH00949303},
\cite[chap.4-6]{zbMATH02113718}. Denote  
\begin{equation}
    \bfB_\Bd = \mrm{diag}(\bfB_{1,\Bd},\dots, \bfB_{\mJ,\Bd})
    \qquad\text{with}\qquad
    \bfB_{j,\Bd}(\bv) = (v_e)_{e\in \Edge_j\setminus\Skel_j}\;\text{for}\;\bv = (v_e)_{e\in\Edge_j}.
\end{equation}
With this notation we have $\bfB_\Bd^\top\bfB_\Bd + \bfB^\top\bfB = \bfId$
which offers a decomposition of unknown vectors into
the degrees of freedom associated to the extended skeleton (labelled "$\Skel$")
and those associated to the interior (labelled "$\In$").
The matrix of the global problem can
then be decomposed accordingly: up to a reordering, it writes as follows
\begin{equation}
  \begin{aligned}
    \bfA =
    \left\lbrack\begin{array}{ll}
    \bfA_{\In \In} &     \bfA_{\In \Bd}\\
    \bfA_{\Bd \In} &     \bfA_{\Bd \Bd}
    \end{array}\right\rbrack
    \quad\text{with}\quad 
    &\bfA_{\In \In}:= \bfB_{\Bd}\bfA\bfB_{\Bd}^\top,
    \quad \bfA_{\In \Bd}:= \bfB_{\Bd}\bfA\bfB^\top,\\
    &\bfA_{\Bd \In}:= \bfB\bfA\bfB_{\Bd}^\top,
    \quad \bfA_{\Bd \Bd}:= \bfB\bfA\bfB^\top.
  \end{aligned}
\end{equation}
The matrix $\bfA_{\Bd \Bd} -
\bfA_{\Bd \In}\bfA_{\In \In}^{-1}\bfA_{\In \Bd}\in \CC^{\#\Skel_{\oplus}\times \#\Skel_{\oplus}}$
is customarily referred to as the Schur complement of $\bfA$ (with respect to
skeleton unknowns). With the help of the Schur complement, the expression
of the scattering matrix becomes simple. 
\begin{lem}\quad\\
  Assume that the matrix $\bfA_{\In\In}$ is invertible, and denote
   $\tilde{\bfA} := \bfT^{-1}(
  \bfA_{\Bd \Bd} - \bfA_{\Bd \In}\bfA_{\In \In}^{-1}
  \bfA_{\In \Bd})$. Then the scattering matrix admits the expression 
  \begin{equation*}
    \bfS = (\tilde{\bfA}-\imath\bfId)^{-1}(\tilde{\bfA}+\imath\bfId).
  \end{equation*}
\end{lem}
\noindent \textbf{Proof:}

Starting from Expression \eqref{FirstExpressionScattering},
pick an arbitrary $\bp\in \Vh(\Skel_\oplus)$ and let us compute the expression
of $\bq = \bfS(\bp)$. Denote $\bv = (\bfA-\imath\bfB^\top\bfT\bfB)^{-1}(\bfA+\imath\bfB^\top\bfT\bfB)
\bfB^\top\bp$ so that $\bq = \bfB(\bv)$. Decomposing into
interior and boundary contributions, with $\bv^\top = \lbr \bv_\In^\top, \bv_\Bd^\top\rbr$,
we have $\bq = \bv_\Bd$ and the linear system
\begin{equation*}
  \begin{aligned}
    & \bfA_{\In\In} \bv_\In + \bfA_{\In\Bd}\bq = \bfA_{\In\Bd}\bp\\
    & \bfA_{\Bd\In} \bv_\In + (\bfA_{\Bd\Bd}-\imath \bfT)\bq = (\bfA_{\Bd\Bd}+\imath\bfT)\bp.
  \end{aligned}
\end{equation*}
Now eliminating the interior unknowns $\bv_\In$ by "Schur complementing" this system
then leads to the identity $(\bfA_\star-\imath\bfT)\bq = (\bfA_\star+\imath\bfT)\bp$
with $\bfA_\star:= \bfA_{\Bd \Bd} - \bfA_{\Bd \In}\bfA_{\In \In}^{-1}
\bfA_{\In \Bd}$. There only remains to multiply on the left by $\bfT^{-1}$
which leads to the expression we were looking for.  \hfill $\Box$
  
\quad\\
The matrix $\bfA_{\In\In}$ is not guaranteed to be invertible. A non-trivial kernel
corresponds to a resonance phenomenon in a local subproblem. This however cannot
occur if the maximum diameter of subdomains is small enough, see e.g. Lemma 11.4 in
\cite{zbMATH02113718}.

The previous lemma delivers the instructive insight that, under appropriate circumstances ($\bfA_{\In\In}$ invertible),
the scattering matrix takes the form of a Cayley
transform. Let us underline however that, even when $\bfA_{\In\In}$ is not invertible,
$\bfA-\imath\bfB^\top\bfT\bfB$ is invertible and the scattering matrix given by
\eqref{FirstExpressionScattering} is properly defined.

\begin{lem}\label{EnergyConservation}\quad\\
  For any $\bp\in \Vh(\Skel_\oplus)$, we have the estimate
  $\Vert \bfS(\bp)\Vert_\bfT\leq \Vert \bp\Vert_\bfT$.
  More precisely,
  recalling the definition of the energy dissipation
  functional~\eqref{DefinitionDissipatedEnergy}, the
  following energy conservation identity holds
  \begin{equation*}
    \Vert \bfS(\bp)\Vert_\bfT^2 + 4\mathcal{P}(\bv) = \Vert \bp\Vert_\bfT^2
    \;\; \text{with} \;\;
    \bv = (\bfA-\imath\bfB^\top\bfT\bfB)^{-1}\bfB^\top\bfT(\bp).
  \end{equation*}
\end{lem}
\noindent \textbf{Proof:}

According to \eqref{DefScatteringMatrix}, we have $\bfS(\bp) = \bp+2\imath \bfB(\bv)$.
Using this expression we have $\Vert \bfS(\bp)\Vert_\bfT^2 =
\Vert \bp+2\imath \bfB\bv\Vert_\bfT^2 = \Vert \bp\Vert_\bfT^2 + 4\Vert \bfB\bv\Vert_\bfT^2
-4\Re e\{\imath (\bp,\bfB \bv)_\bfT\}$. On the other hand, the very definition of
$\bv$ directly yields $-\overline{\bv}^\top\bfA\bv+\imath (\bfB\overline{\bv})^\top\bfT(\bfB\bv) =
-(\bfB\overline{\bv})^\top\bfT(\bp)$  
which rewrites $\mathcal{P}(\bv) := -\Im m\{\overline{\bv}^\top\bfA\bv\} =
-\Vert \bfB\bv\Vert_\bfT^2+\Re e\{\imath (\bp,\bfB\bv)_\bfT \}$. From this follows the
desired energy conservation identity and, since $\mathcal{P}(\bv)\geq 0$ according
to \eqref{DefinitionDissipatedEnergy}, we also deduce
$\Vert \bfS(\bp)\Vert_\bfT\leq \Vert \bp\Vert_\bfT$. \hfill $\Box$

\quad\\
From the previous identity, we deduce that $\Vert \bfS(\bp)\Vert_\bfT\leq \Vert \bp\Vert_\bfT$
for all $\bp\in\Vh(\Skel_\oplus)$ i.e. the scattering matrix is non-expansive.
The previous energy conservation result actually paves the way to proving
the invertibility of the matrix of \eqref{SkeletonFormulation}.

\begin{prop}\label{prop:invertibility}\quad\\
  The matrix $\bfId+\bfPi\bfS\in \CC^{\#\Skel_{\oplus}\times
  \#\Skel_{\oplus}}$ is invertible.
\end{prop}
\noindent \textbf{Proof:}

We need to show that $\ker(\bfId+\bfPi\bfS) = \{0\}$.
Pick any $\bp\in \Vh(\Skel_\oplus)$ satisfying
$(\bfId+\bfPi\bfS)\bp = 0$
and set $\bu = (\bfA-\imath\bfB^\top\bfT\bfB)^{-1}\bfB^\top\bfT(\bp)$. As we already
mentioned before, we have $\bfS(\bp) = \bp+2\imath\bfB(\bu)$ so the following
equations hold
\begin{equation*}
    (\bfA-\imath\bfB^\top\bfT\bfB)\bu = \bfB^\top\bfT(\bp),
    \qquad\text{and}\qquad
    \bp+\bfPi(\bp+2\imath\bfB\bu) = 0.
\end{equation*}
This means that the pair $(\bu,\bp)$ must be solution to \eqref{Formulation4}
with $\bff = 0$.
The latter problem was shown to be equivalent to \eqref{Formulation1}.
This implies that $\bu\in \Vh(\Edge_\oplus)$ solves~\eqref{Formulation1} with
$\bff = 0$, and that $\bu_\Omega = (\bfR^\top\bfR)^{-1}\bfR^\top\bu$
solves~\eqref{Formulation0} with \(\bff_{\Omega}=0\).
Since~\eqref{Formulation0} was assumed uniquely solvable,
we deduce that $\bu_\Omega = 0\Rightarrow \bu = \bfR(\bu_\Omega) = 0$.
This implies that $\bfB^\top\bfT\bp = 0$.
Since $\bfB\bfB^\top = \bfId$, multiplying on the left by 
$\bfB$ yields $\bfT\bp = 0$ hence $\bp = 0$ as $\bfT$ is assumed symmetric
positive definite. This proves the invertibility of $\bfId+\bfPi\bfS$.
\hfill $\Box$

\begin{prop}\label{prop:well-posedness}\quad\\
  All eigenvalues of $\bfId+\bfPi\bfS$ belong to the punctured
  disk $\{\lambda\in \CC\setminus\{0\},\;\vert 1-\lambda\vert\leq 1 \}$. Moreover
  we have $\Vert (\bfId+\bfPi\bfS)\bp\Vert_\bfT\leq 2\Vert \bp\Vert_\bfT$ for
  all $\bp\in\Vh(\Skel_\oplus)$, and there exists a constant $\alpha>0$ such
  that,
  for all $\bp\in \Vh(\Skel_\oplus)$,
  \begin{equation*}
    \Re e\{(\bp,(\bfId+\bfPi\bfS)\bp)_\bfT\}\geq \alpha \Vert \bp\Vert_\bfT^2.
  \end{equation*}
\end{prop}
\noindent \textbf{Proof:}

The property on the location of eigenvalues and the upper bound stem directly from
the inequality $\Vert \bfPi\bfS(\bp)\Vert_{\bfT}\leq \Vert \bp\Vert_{\bfT}$
(see Lemma~\ref{PropertiesExchangeMatrix} and Lemma~\ref{EnergyConservation})
as well as the invertibility of $\bfId+\bfPi\bfS$
from Proposition~\ref{prop:invertibility}.
Next, set
\begin{equation}
  \alpha  = \inf_{\bq\in \Vh(\Skel_\oplus)\setminus\{0\}}  \frac{\Re e\{(\bq,(\bfId+\bfPi\bfS)\bq)_\bfT\}}{
  \Vert \bq\Vert_\bfT^2}.
\end{equation}
Take $\bp\in\Vh(\Skel_\oplus)\setminus\{0\}$ with $\Vert \bp\Vert_\bfT = 1$ and
$\Re e\{(\bp,(\bfId+\bfPi\bfS)\bp)_\bfT\} = \alpha$. Applying Cauchy-Schwarz inequality
together with
Lemma~\ref{PropertiesExchangeMatrix} and Lemma~\ref{EnergyConservation}
already gives $\alpha  = 1 + \Re e\{(\bp,\bfPi\bfS\bp)_\bfT\}\geq 1- \Vert
\bp\Vert_\bfT\Vert \bfPi\bfS(\bp)\Vert_\bfT = 1- \Vert
\bfS(\bp)\Vert_\bfT\geq 1- \Vert \bp\Vert_\bfT^2 = 0$ i.e. $\alpha \geq 0$.
Next, proceed by contradiction, assuming $\Re e\{(\bp,(\bfId+\bfPi\bfS)\bp)_\bfT\} = \alpha = 0$.
According to
Lemma~\ref{PropertiesExchangeMatrix} and Lemma~\ref{EnergyConservation}
again, we have
\begin{equation}
  \begin{aligned}
    \Vert \bp\Vert_\bfT^2\geq
    & \Vert \bfS(\bp)\Vert_\bfT^2 = \Vert \bfPi\bfS(\bp)\Vert_\bfT^2
     =  \Vert (\bfId+\bfPi\bfS)\bp - \bp\Vert_\bfT^2 \\
     &  =  \Vert (\bfId+\bfPi\bfS)\bp\Vert_\bfT^2 +\Vert \bp\Vert_\bfT^2
     - 2\Re e\{(\bp,(\bfId+\bfPi\bfS)\bp)_\bfT\}
     =  \Vert (\bfId+\bfPi\bfS)\bp\Vert_\bfT^2 +\Vert \bp\Vert_\bfT^2.
  \end{aligned}
\end{equation}
From this we finally conclude that $\Vert (\bfId+\bfPi\bfS)\bp\Vert_\bfT = 0$ which
shows that $\bp = 0$ since $\bfId+\bfPi\bfS$ was proved invertible. This contradicts
$\Vert \bp\Vert_\bfT = 1$ so we finally conclude that $\alpha >0$
necessarily. \hfill $\Box$

\quad\\
The previous result directly implies the convergence of standard fixed point
algorithms such as the damped Richardson algorithm, hence necessarily the
convergence of the restarted \textsc{gmres} solver.

\section{Concrete definitions of transmission matrices}\label{sec:transmission}

In the present section we examine and discuss two concrete choices of transmission matrices. 

\subsection{Zeroth-order transmission matrices}
We discuss a first choice of transmission matrix based on the $\mathrm{L}^2$ scalar
product of tangential traces.
We assume here that
\(\Skel_{j}\) contains only edges of \(\partial \Omega_{j}\).
In Section~\ref{sec:num_pie} we used 
\(\Skel_{j} = \EdgeStar \cap \Edge_{j} \)
while in the rest of the numerical experiments 
\(\Skel_{j}\) contains all edges of \(\partial \Omega_{j}\).
The transmission matrix
$\bfT = \mrm{diag}(\bfT_1,\dots,\bfT_\mJ)$ where the entries of the matrices $\bfT_j:\Vh(\Gamma_j)\to\Vh(\Gamma_j)$
are given by
\begin{equation}\label{PurelyLocalImpedance}
  (\bfT_j)_{e,f} = \int_{\Gamma_j}(\kappa/\check{\eta}_{j})(\bvphi_e\times\bn_j)\cdot(\bvphi_f\times\bn_j) d\sigma.
\end{equation}
The function $\check{\eta}_{j}:\Gamma_{j}\to (0,+\infty)$ can be
chosen arbitrarily.
With such a choice of transmission matrix, local problems amount to numerically
solving Maxwell problems in each subdomain $\Omega_j$ with the first order
absorbing boundary condition
$\bn_j\times \bE\times \bn_j - \check{\eta}_{j}\,\bH\times\bn_j=\bg$
on $\Gamma_j$,
 for some $\bg$.
A common choice for \(\check{\eta}_{j}\) is to take the value of
$\Re e\{\sqrt{\mu/\epsilon}\}$.
This quantity might be discontinuous across \(\Gamma_{j}\), when the
coefficients \(\epsilon\) and \(\mu\) are non-constant in the domain \(\Omega\).
In this case, an average over neighboring mesh cells is commonly performed
to get a single value at the interface cell.
We point out importantly that in our approach, this is \emph{not} a requirement.
As a result, our approach provides much more flexibility and can handle
discontinuities seamlessly.

This choice of transmission matrix corresponds to the strategy originally used
in the work of Despr\'es \cite{zbMATH00863010}, assuming that
\(\check{\eta}_{j}\) takes the same value from each side of \(\Gamma_{j}\)
so that~Hypothesis~\eqref{BlockDiagonalForm} is satisfied.
This work and the variants considered so far in the literature discards the
issue raised by the presence of cross-points by adopting a different
discretization scheme: a mixed hybrid discretization \cite{zbMATH00440677}
where the degrees of freedom are associated to the faces of each tetrahedron
and can thus be easily exchanged by a simple swap.
On the contrary, the approach we adopt here is able to deal with the presence
of degrees of freedom at cross-points even for N\'ed\'elec finite elements. To
be more specific, for onion skin domain decompositions, we only have to handle
single interfaces with edges of multiplicity two and Choice
\eqref{PurelyLocalImpedance} fits the situation described at the end of \S
\ref{ExplicitExpression}. Consequently, the communication matrix is given
explicitly by Formula \eqref{ExpressionExplicitSwapping} and our method is a
simple extension of Despr\'es' method to any conformal finite element method.
However, by introducing a more general communication matrix, our theory
allows to deal with domain decomposition with simple transmission quantities
involving degrees of freedom of multiplicity more than two for N\'ed\'elec's
elements. This appears to be new.

\subsection{Schur complement based transmission matrix}\label{sec:impedance_Schur}
We shall now examine an alternative possible choice of transmission matrix based on the
Schur complement associated to the solution of some auxiliary strongly coercive
problem. We dedicate a whole section to this particular transmission matrix because it
appears as one of the most efficient choices.

We first need to consider a subset $\Omega'\subset \Omega$ obtained as union of
a subset of elements of the triangulation $\Tri(\Omega')\subset \Tri(\Omega)$
and such that $\Omega' = \cup_{\tau\in \Tri(\Omega')}\overline{\tau}$.
Setting $\Omega_j':=\Omega'\cap \Omega_j$,
we have $\overline{\Omega}\!\,' = \cup_{j=1}^{\mJ}\overline{\Omega}\!\,_{j}' $.
Next denote $\Edge'$ the collection of edges of
$\Tri(\Omega')$, as well as $\Edge'_j$ those belonging to $\Tri(\Omega'_j)$.
We make the following important assumption that
$\Tri(\Omega')$ is selected so as to guarantee that $\Gamma\subset \Edge'$ and
$\Gamma_j\subset \Edge'_j$.
The subset $\Omega'$ will be the computational domain for our auxiliary
problem. It shall typically consist in layers of elements
surrounding the skeleton \eqref{DefThinSkeleton} of the subdomain
decomposition, see Figure~\ref{FigSkeleton}(b).
In each subdomain $\Omega'_j$ we consider a bilinear form
\begin{equation}\label{eq:schur_complement_bilinear_form}
  \begin{aligned}    
    & c_j(\bu,\bv) :=\int_{\Omega_j'}\Re e\{\mu^{-1}\}\bfcurl(\bu)\cdot\bfcurl(\overline{\bv})
    + \kappa^2\Re e\{\epsilon\} \bu\cdot\overline{\bv} \;d\bx\\
    & \textcolor{white}{c_j(\bu,\bv)} + \int_{\partial\Omega_j'\setminus(\partial\Omega_j\setminus\partial\Omega)}\Re e\{\kappa/\eta\}(\bu\times\bn_j')
    \cdot(\overline{\bv}\times \bn_j')\;d\sigma\\    
    & \bfC_j:\Vh(\Edge'_j)\to \Vh(\Edge'_j), \quad (\bfC_j)_{e,f} := c_j(\bvphi_f,\bvphi_e).
  \end{aligned}
\end{equation}
where $\bn_j'$ refers to the vector field normal to $\partial\Omega_j'$. We also set
$\bfC:=\mrm{diag}(\bfC_1,\dots,\bfC_\mJ)$. By construction this is a symmetric positive
definite matrix. Next we separate unknowns located on
the skeleton from other unknowns by means of restriction matrices and we define the
auxiliary matrix $\bfB':\Vh(\Edge'_\oplus)\to \Vh(\Gamma_\oplus)$ by
\begin{equation}
   \bfB' := \mrm{diag}(\bfB_1',\dots, \bfB_\mJ')
   \qquad\text{with}\qquad
   \bfB_j'(\bv) :=  (v_e)_{e\in \Gamma_j}\quad \bv=(v_e)_{e\in \Edge'_j}\in \Vh(\Edge'_j).
\end{equation}
The transmission matrix that we propose to consider here is the Schur complement associated to the elimination of interior
unknowns in the matrix $\bfC$ defined above. To be more specific we consider the matrix $\bfT:\Vh(\Gamma_\oplus)\to \Vh(\Gamma_\oplus)$
defined by
\begin{equation}\label{ExpressionSchurComplement}
    \bfT(\bu_{\Bd}) = \bq\quad \text{where} \quad
    \; (\bv,\bq)\in \Vh(\Edge'_\oplus)\times \Vh(\Gamma_\oplus)
    \;\;\text{solves}\;\;
    \left\lbr
    \begin{array}{cc}
      \bfC  & -(\bfB')^{\top}\\
      \bfB' & 0
    \end{array}
    \right\rbr\cdot
    \left\lbr
    \begin{array}{l}
      \bv\\\bq
    \end{array}
    \right\rbr =
    \left\lbr
    \begin{array}{l}
      0\\\bu_{\Bd}
    \end{array}
    \right\rbr.
\end{equation}
As a Schur complement of a SPD matrix, it is itself SPD and is thus a valid candidate for the construction presented
in Section~\ref{sec:OrtogonalProjection}. To obtain an expression for the final system to be considered in the global DDM
strategy, we need to combine \eqref{ExpressionSchurComplement} with \eqref{Formulation4}. In this process, one has
$\bu_{\Bd} = -\imath\bfB\bu - \bp$ which leads to the system
\begin{equation}\label{ExpressionSchurComplement2}
  \begin{aligned}
    & \text{Find}\quad (\bu,\bv)\in \Vh(\Edge_\oplus)\times\Vh(\Edge_\oplus'), \;\;
    (\bp,\bq)\in \Vh(\Gamma_\oplus)\times\Vh(\Gamma_\oplus)
    \;\text{such that}\\
    &
    \left\lbr\begin{array}{ccc}
    \bfA & 0 & \;\bfB^\top\\
    0    & -\imath\bfC & -\bfB'^\top\\
    \bfB & -\bfB'      & 0
    \end{array}\right\rbr\cdot
    \left\lbr\begin{array}{c}
    \bu \\ \bv \\ \bq
    \end{array}\right\rbr =
    \left\lbr\begin{array}{c}
    \bff \\ 0 \\ \imath\bp
    \end{array}\right\rbr,
    \quad \text{and}\quad \bp = - \boldsymbol{\Pi}(\bp+2\imath\bfB\bu).
  \end{aligned}\\ 
\end{equation}
Of course we also have to discuss actual computation  of the matrix $\boldsymbol{\Pi} = 2\bfP-\bfId$ since,
according to \eqref{ExchangeMatrixDef} and  \eqref{NeumannNeumannStrategy}, it involves matrix-vector
products for both $\bfT$ and $\bfT^{-1}$. Matrix-vector product by  $\bfT$ can be treated based on
\eqref{ExpressionSchurComplement}. Matrix-vector by $\bfT^{-1}$ can be computed
using the identity $\bfT^{-1} = \bfB'\bfC^{-1}\bfB'^{\top}$.

Despite their appearing in the right hand side, $\bu,\bp$ are unknowns of
\eqref{ExpressionSchurComplement2}, and only $\bff$ is a source term.
We arranged a system of equations like in \eqref{ExpressionSchurComplement2}
in the perspective of an iterative solution procedure. In practice, for Schur complement
based transmission matrices as discussed in the present paragraph, the linear system
appearing in the left hand side of \eqref{ExpressionSchurComplement2} is the one to be dealt
with at each iteration for applying the scattering matrix $\bfS$ defined
in \eqref{DefScatteringMatrix}. This can be achieved in parallel thanks to the subdomain-wise
block diagonal structure of the left hand side of \eqref{ExpressionSchurComplement2}.

We advocate the design of transmission matrices like \eqref{ExpressionSchurComplement} because, under
technical assumptions, it is shown \cite{parolin:tel-03118712,claeys2020robust} that the coercivity
constant of $\bfId + \boldsymbol{\Pi}\bfS$ is bounded from below independently of the meshsize
which leads to robust convergence of linear iterative solvers applied to the skeleton formulation
\eqref{SkeletonFormulation}.

\section{Algorithms}\label{sec:algorithms}

We wish now to describe in more concrete terms the practical implementation of
the method.
Our emphasis is on the parallel nature of the algorithms, in particular the
\texttt{for} loops over the \(\mJ\) subdomains are written explicitly and can be
parallelized.
Recall that the problem that is solved in practice is
Problem~\eqref{SkeletonFormulation} which is posed on the extended skeleton.

\subsection{General algorithms}

We first provide the general forms of the algorithms by which we mean the
definitions of the algorithms that can be applied for any generic scalar
product \(\bfT\) given by a family of local contributions \(\bfT_{j}\).
Such procedures are in particular well-adapted to the Despr\'es transmission matrix.

\paragraph{Richardson algorithm.}
The damped Richardson algorithm is first considered, with damping parameter
denoted by \(r\).
Besides the definitions of the restriction matrices, recall in particular the
definitions of the local contributions \(\bfA_{j}\) and \(\bff_{j}\)
in~\eqref{LocalContrib}.
The general form of the Richardson algorithm is then given in
Algorithm~\ref{alg:general_Richardson}.

\begin{algorithm}
  \caption{General form of the Richardson algorithm}\label{alg:general_Richardson}
  \begin{algorithmic}[1]
    \For{$j=1,\dots, \mJ$} \Comment{Initialisation}
    \State%
    $\bp_{j} = 0$ \Comment{size: $\#\Skel_{j}$}
    \State%
    $\bu_{j} = {(\bfA_{j}-\imath \bfB_{j}^{\top}\bfT_{j}\bfB_{j})}^{-1}\bff_{j}$
    \label{lin:general_Richardson_init_local_solve}
    \Comment{Local solve (size: $\#\Edge_{j}$)}
    \EndFor%
    \For{$n=1,\dots,n_{\max}$}
    \State%
    $\bg = 0$ \Comment{size: $\#\Skel$}
    \For{$j=1,\dots, \mJ$}
    \State%
    $\bg = \bg + \bfQ_{j}^{\top}\bfT_{j}(\bp_{j} + 2\imath\bfB_{j}\bu_{j})$
    \Comment{Local scattering}
    \label{lin:general_Richardson_apply_T}
    \EndFor%
    \State%
    $\bv = \left(\bfQ^{\top}\bfT\bfQ\right)^{-1}\bg$
    \Comment{Solved using PCG (size: $\#\Skel$)}
    \label{lin:general_Richardson_pcg}
    \For{$j=1,\dots, \mJ$}
    \State%
    $\bp_{j} = \bp_{j} + 2r (\imath \bfB_{j}\bu_{j} -\bfQ_{j}\bv)$
    \State%
    $\bu_{j} = {(\bfA_{j}-\imath \bfB_{j}^{\top}\bfT_{j}\bfB_{j})}^{-1}{(\bfB_{j}^{\top}\bfT_{j}\bp_{j}+\bff_{j})}$
    \Comment{Local solve (size: $\#\Edge_{j}$)}
    \label{lin:general_Richardson_loop_local_solve}
    \EndFor%
    \EndFor%
  \end{algorithmic}
\end{algorithm}

Of course, in the above algorithm (and in the algorithms below) the inverse
matrices, namely
\({(\bfA_{j}-\imath \bfB_{j}^{\top}\bfT_{j}\bfB_{j})}^{-1}\)
are not actually assembled. Instead, each matrix
\(\bfA_{j}-\imath \bfB_{j}^{\top}\bfT_{j}\bfB_{j}\)
is factorized (offline precomputations) and the inversion of the linear system
is performed in the course of the iterations using forward and backward
substitution.

Besides, as explained above, the projection problem in
Algorithm~\ref{alg:general_Richardson} appearing in
Line~\ref{lin:general_Richardson_pcg}
is performed using a preconditioned CG algorithm.
To define the PCG algorithm, it suffices to provide a definition for a
matrix-vector product routine for the problem matrix
\(\bfQ^{\top}\bfT\bfQ\) as well as the preconditioner matrix \(\bfM\),
see~\eqref{NeumannNeumannStrategy}.
Albeit the fact that such routines are straightforward, the procedure are
respectively provided in Algorithm~\ref{alg:general_matvec_CG} and in
Algorithm~\ref{alg:general_matvec_precond_CG}
to stress in particular that they are fully parallel.
Notice that the matrix \(\bfD\) is diagonal.

\begin{center}
\begin{minipage}{0.48\linewidth}
\begin{algorithm}[H]
  \caption{Matrix-vector product for CG}\label{alg:general_matvec_CG}
  \begin{algorithmic}[1]
    \Input
    $\bg$
    \Comment{size: $\#\Skel$}
    \State%
    $\bq = 0$
    \Comment{size: $\#\Skel$}
    \For{$j=1,\dots, \mJ$}
    \State%
    $\bq = \bq + \bfQ_{j}^{\top}\bfT_{j}\bfQ_{j}\bg$
    \label{lin:general_matvec_CG_apply_T}
    \EndFor%
    \Output
    $\bq$
  \end{algorithmic}
\end{algorithm}
\end{minipage}
\hfill
\begin{minipage}{0.47\linewidth}
\begin{algorithm}[H]
  \caption{CG preconditioner}\label{alg:general_matvec_precond_CG}
  \begin{algorithmic}[1]
    \Input
    $\bq$
    \Comment{size: $\#\Skel$}
    \State%
    $\bq = \bfD\bq$
    \State%
    $\bp = 0$
    \Comment{size: $\#\Skel$}
    \For{$j=1,\dots, \mJ$}
    \State%
    $\bp = \bp + \bfQ_{j}^{\top}\bfT_{j}^{-1}\bfQ_{j}\bq$
    \label{lin:general_matvec_precond_CG_apply_T}
    \EndFor%
    \State%
    $\bp = \bfD\bp$
    \Output
    $\bp$
  \end{algorithmic}
\end{algorithm}
\end{minipage}
\end{center}

\paragraph{\texttt{GMRES} algorithm.}
The Richardson algorithm is rarely used in practice and Krylov methods are the
preferred choice in real-life applications.
Since the wave propagation problems yields non-symmetric problems, one will
typically resort to the \textsc{gmres} algorithm.

To define the \textsc{gmres} algorithm, it suffices to provide a definition for a
right-hand side and a matrix-vector product routine.
The right-hand side is denoted by \(\bb\) (see~\eqref{DefScatteringMatrix})
and can be computed (offline)
according to Algorithm~\ref{alg:general_RHS}.
The matrix-vector product procedure, which
takes as input a vector \(\bp\) and outputs a vector \(\bq\), is given in
Algorithm~\ref{alg:general_matvec}.
Again, the projection problem
appearing in Line~\ref{lin:general_gmres_rhs_pcg}
of Algorithm~\ref{alg:general_RHS}
and in Line~\ref{lin:general_gmres_pcg}
of Algorithm~\ref{alg:general_matvec}
is performed using the same preconditioned CG algorithm that was defined for
the Richardson algorithm.

\begin{algorithm}
  \caption{RHS computation for \textsc{gmres}}\label{alg:general_RHS}
  \begin{algorithmic}[1]
    \State%
    $\bb = 0$ \Comment{size: $\#\Skel_{\oplus}}$
    \State%
    $\bg = 0$ \Comment{size: $\#\Skel$}
    \For{$j=1,\dots, \mJ$}
    \State%
    $\bu_{j} = {(\bfA_{j}-\imath \bfB_{j}^{\top}\bfT_{j}\bfB_{j})}^{-1}\bff_{j}$
    \Comment{Local solve (size: $\#\Edge_{j}$)}
    \label{lin:general_RHS_local_solve}
    \State%
    $\bb = \bb + 2\imath\, \bfQ\bfQ_{j}^{\top}\bfB_{j}\bu_{j}$
    \Comment{size: $\#\Skel_{j}$}
    \State%
    $\bg = \bg + 2\imath\, \bfQ_{j}^{\top}\bfT_{j}\bfB_{j}\bu_{j}$
    \label{lin:general_RHS_apply_T}
    \EndFor%
    \State%
    $\bv = \left(\bfQ^{\top}\bfT\bfQ\right)^{-1}\bg$
    \Comment{Solved using PCG (size: $\#\Skel$)}
    \label{lin:general_gmres_rhs_pcg}
    \State%
    $\bb = \bb - 2\, \bfQ\bv$
    \Output
    $\bb$
  \end{algorithmic}
\end{algorithm}

\begin{algorithm}
  \caption{Matrix-vector product for \textsc{gmres}}\label{alg:general_matvec}
  \begin{algorithmic}[1]
    \Input
    $\bp$
    \State%
    $\bq = 0$ \Comment{size: $\#\Skel_{\oplus}}$
    \State%
    $\bg = 0$ \Comment{size: $\#\Skel$}
    \For{$j=1,\dots, \mJ$}
    \State%
    $\bu_{j} = {(\bfA_{j}-\imath \bfB_{j}^{\top}\bfT_{j}\bfB_{j})}^{-1}{(\bfB_{j}^{\top}\bfT_{j}\bp_{j})}$
    \Comment{Local solve (size: $\#\Edge_{j}$)}
    \label{lin:general_matvec_local_solve}
    \State%
    $\bq = \bq - 2\imath\, \bfQ\bfQ_{j}^{\top}\bfB_{j}\bu_{j}$
    \Comment{size: $\#\Skel_{j}$}
    \State%
    $\bg = \bg + \bfQ_{j}^{\top}\bfT_{j}(\bp_{j}+2\imath\bfB_{j}\bu_{j})$
    \label{lin:general_matvec_apply_T}
    \Comment{Local scattering}
    \EndFor%
    \State%
    $\bv = \left(\bfQ^{\top}\bfT\bfQ\right)^{-1}\bg$
    \Comment{Solved using PCG (size: $\#\Skel$)}
    \label{lin:general_gmres_pcg}
    \State%
    $\bq = \bq + 2\, \bfQ\bv$
    \Output
    $\bq$
  \end{algorithmic}
\end{algorithm}

\subsection{Algorithms with the Schur complement based transmission matrix}

We now turn to the particular case where one uses a Schur complement based
transmission matrix and explain how the above algorithms need to be
modified.
As we already explained, the algorithms can be written so that no dense
matrix is involved (i.e. the Schur complement is not performed in practice),
albeit the underlying non-local nature of the transmission operator.
This is particularly important for efficiency considerations because otherwise
the naive implementation of the method requires the computation and storage of
dense matrices as well as the solution to hybrid sparse-dense linear systems
for which many factorization routines may struggle.

Before describing the algorithms let us define
\begin{equation}
  \tilde{\bfC}_{j} := 
  \left\lbr
  \begin{array}{cc}
    \bfC_{j}  & -(\bfB_{j}')^{\top}\\
    \bfB_{j}' & 0
  \end{array}
  \right\rbr,
  \qquad\text{and}\qquad
  \tilde{\bfA}_{j} := 
  \left\lbr
  \begin{array}{ccc}
    \bfA_{j} & 0 & \;\bfB_{j}^\top\\
    0    & -\imath\bfC_{j} & -\bfB_{j}'^\top\\
    \bfB_{j} & -\bfB_{j}'      & 0
  \end{array}
  \right\rbr.
\end{equation}
The matrices \(\bfC_{j}\), \(\tilde{\bfC}_{j}\) and \(\tilde{\bfA}_{j}\)
are fully sparse matrices than can be factorized (offline).
In the algorithms their inverses will appear, which correspond
in practice to forward and backward substitutions.
The matrix \(\tilde{\bfC}_{j}\)
has size $\#\Edge_{j}'+\#\Skel_{j}$
and \(\tilde{\bfA}_{j}\) has size
$\#\Edge_{j}+\#\Edge_{j}'+\#\Skel_{j}$.

\paragraph{Richardson algorithm.}

We now give the modifications regarding the Richardson algorithm for the Schur
complement based transmission matrix.
The local solve appearing in
line~\ref{lin:general_Richardson_init_local_solve}
of Algorithm~\ref{alg:general_Richardson}
is replaced by
$(\bu_{j},\bv_{j},\bq_{j})^{\top} = \tilde{\bfA}_{j}^{-1}(\bff_{j},0,0)^{\top}$
and the one of line~\ref{lin:general_Richardson_loop_local_solve}
by
$(\bu_{j},\bv_{j},\bq_{j})^{\top} = \tilde{\bfA}_{j}^{-1}(\bff_{j},0,\imath\bp_{j})^{\top}$.
The computation of the quantity
$\bfT_{j}(\bp_{j} + 2\imath\bfB_{j}\bu_{j})$
in line~\ref{lin:general_Richardson_apply_T}
is replaced by the quantity $\bq_{j}$
computed as
$(\bv_{j},\bq_{j})^{\top} = \tilde{\bfC}_{j}^{-1}(0,\bp_{j} + 2\imath\bfB_{j}\bu_{j})^{\top}$.

Again in this particular case, the projection problem
is solved using a preconditioned CG algorithm.
The computation of the quantity
$\bfT_{j}\bfQ_{j}\bg$
in line~\ref{lin:general_matvec_CG_apply_T}
of Algorithm~\ref{alg:general_matvec_CG}
is replaced by the quantity $\bq_{j}$
computed as
$(\bv_{j},\bq_{j})^{\top} = \tilde{\bfC}_{j}^{-1}(0,\bfQ_{j}\bg)^{\top}$.
The computation of the quantity
$\bfT_{j}^{-1}\bfQ_{j}\bg$
in line~\ref{lin:general_matvec_precond_CG_apply_T}
of Algorithm~\ref{alg:general_matvec_precond_CG}
is replaced by the quantity
$\bfB_{j}'\bfC_{j}^{-1}\bfB_{j}'^{\top}\bfQ_{j}\bq$.

\paragraph{\texttt{GMRES} algorithm.}

We now give the modifications regarding the Krylov algorithm for the Schur
complement based transmission matrix.
The local solve appearing in
line~\ref{lin:general_RHS_local_solve}
of Algorithm~\ref{alg:general_RHS}
is replaced by
$(\bu_{j},\bv_{j},\bq_{j})^{\top} = \tilde{\bfA}_{j}^{-1}(\bff_{j},0,0)^{\top}$
and the one of line~\ref{lin:general_matvec_local_solve}
of Algorithm~\ref{alg:general_matvec}
by
$(\bu_{j},\bv_{j},\bq_{j})^{\top} = \tilde{\bfA}_{j}^{-1}(0,0,\imath\bp_{j})^{\top}$
The computation of the quantity
$\bfT_{j}(\bp_{j} + 2\imath\bfB_{j}\bu_{j})$
in line~\ref{lin:general_RHS_apply_T}
of Algorithm~\ref{alg:general_RHS}
is replaced by the quantity $\bq_{j}$
computed as
$(\bv_{j},\bq_{j})^{\top} = \tilde{\bfC}_{j}^{-1}(0,\bfB_{j}\bu_{j})^{\top}$
and the one in line~\ref{lin:general_matvec_apply_T}
of Algorithm~\ref{alg:general_matvec}
is replaced by the quantity $\bq_{j}$
computed as
$(\bv_{j},\bq_{j})^{\top} = \tilde{\bfC}_{j}^{-1}(0,\bp_{j} + 2\imath\bfB_{j}\bu_{j})^{\top}$.

\section{Numerical experiments}\label{sec:numerics}

We present now a sequence of numerical experiments supporting the previous analysis
and illustrating a few features of the novel approach.

In all our test cases we consider a transmission problem either in a disk (in
2D) or in a ball (in 3D).
We set \(\bJ\equiv 0\) and consider a source that comes from an inhomogeneous
condition on the exterior boundary
\(\eta \bJ_\sigma =  \bn\times [\bE^{\mathrm{inc}}\times \bn] - \eta\,\bH^{\mathrm{inc}}\times\bn\)
where 
\((\bE^{\mathrm{inc}},\bH^{\mathrm{inc}})\)
corresponds to an incoming
plane wave i.e.\
\(\bE^{\mathrm{inc}} = \mathbf{x} \mapsto
\hat{\mathbf{y}}e^{\imath\kappa \mathbf{x} \cdot \hat{\mathbf{x}}}\)
with \((\hat{\mathbf{x}}, \hat{\mathbf{y}})\) the unit vectors in cartesian
coordinates.
The propagation medium is always considered homogeneous with coefficients
\(\mu_r\equiv\epsilon_r\equiv\eta_r\equiv 1\), except in
Section~\ref{sec:num_heterogeneity} where we consider a medium with varying
coefficients \(\epsilon_r\) and \(\mu_r\).
We consider the two transmission matrices that were described in
Section~\ref{sec:transmission} (except in Section~\ref{sec:num_pie} where an
alternative to the Schur complement approach is also considered).
In Section~\ref{sec:num_pie}, the extended skeleton is chosen to be
\(\Skel=\EdgeStar\) (see Figure~\ref{FigSkeleton}(a))
while in the subsequent numerical tests it also includes edges of multiplicity one 
on the physical boundary \(\partial\Omega\) (see Figure~\ref{FigSkeleton}(c)).
While this is not a requirement, for the simplicity of the implementation, the
Schur complement based transmission matrix is constructed in most of our
experiments with \(\Omega'=\Omega\).
The only exception to this rule is the results given at the end of
Section~\ref{sec:num_pie} where we specifically studied an alternative, namely a
much smaller region in the vicinity of the interfaces.

We will present some results where Problem~\eqref{SkeletonFormulation}
is solved using either a damped Richardson iteration scheme (with relaxation
parameter \(r=1/2\)) or a restarted \textsc{gmres} algorithm (with a restart every
\(20\) iterations except in Section~\ref{sec:num_pie} where it is
every \(5\) iterations).
All numerical errors reported
(including the relative error represented in convergence histories) 
are computed between the exact discrete solution of the original (undecomposed)
linear system and the volume broken solutions computed at each iteration
of the iterative solvers. The norm used is the \(\kappa\)-dependent energy norm
which corresponds to the following \(\bfH(\bfcurl)\) norm
\(\|\cdot\|^2 :=
\|\cdot\|^2_{\mathbf{L}^{2}} + \kappa^{-2}\|\bfcurl \cdot\|^2_{\mathbf{L}^{2}}\).

The numerical results were obtained using in-house demonstration codes built
to test the approach.
Meshes were obtained using \textsc{Gmsh}~\cite{Geuzaine2009} and (unless
specified otherwise) partitioned using the automatic graph partitioner
\textsc{Metis}~\cite{Karypis1999}.
The code is mainly sequential (the inherent parallel nature of
the algorithm is not exploited) and is of proof-of-concept nature.
For these reasons, no run times will be reported and we compare different
methods with respect to iteration counts only.
One shall bear in mind though that the cost per iteration is different for each
method.

\subsection{Pie-like splitting}\label{sec:num_pie}

We propose a first test case which aims at illustrating the interest of the
proposed approach.  Our purpose is to give evidence that straightforward
generalizations of more standard methods proposed in the literature, in
particular~\cite{zbMATH01488462,zbMATH07197799}, are not adequate in presence
of cross-points, even in the case where no degrees of freedom are attached to
the cross-points.

In this test case, the unit disk in 2D is regularly split (this is a
geometrically based partitioning, not using an automatic graph partitioner)
into \(\mJ\) pie wedges pointing at the center of the disk.  Therefore, by
construction, there are \(\mJ\) boundary cross-points and one single interior
cross-point (the center of the disk) which is shared by all subdomains.
Since we use (low order) N\'ed\'elec edge finite elements, no degrees of
freedom are attached to the cross-points.
Yet, the numerical results of this section will highlight that,
already in this seemingly simple setting, robustness and uniform convergence
with respect to the discretization parameter in the presence of this interior
cross-point can only be tackled by using a transmission matrix
\(\bfT\) representing a non-local operator
 together with the associated non-local communication matrix \(\bfPi\).

As already mentioned, we consider in this section an alternative to the Schur
complement approach, in addition to the two transmission matrices that were
described in Section~\ref{sec:transmission}.
The difference lies in the location of the degrees of freedom against which the
Schur complement is performed.
In the approach of Section~\ref{sec:impedance_Schur}, they are considered in
the full subdomain boundary. Here we consider also the case where the Schur
complement is performed against each interface (between two subdomains)
independently.
This equivalently amounts to setting to zero off-diagonal blocks
that couple two distinct interfaces in the matrix \(\bfT\) defined in
Section~\ref{sec:impedance_Schur}.
The end result is a block diagonal matrix \(\bfT\) with the number of
blocks corresponding to the number of subdomains in the case of the
matrix of Section~\ref{sec:impedance_Schur}, and to
the number of interfaces in the alternative case considered in addition here.
In particular, this interfaced-based non-local \(\bfT\)
fits the situation described at the end of Section~\ref{ExplicitExpression} and
the communication matrix is given explicitly by
Formula~\eqref{ExpressionExplicitSwapping}.
In some sense, the use of this matrix is the most straightforward extension of
already established approaches akin to~\cite{zbMATH01488462,zbMATH07197799}.
We included this transmission matrix in the numerical results to provide
numerical evidence that the matrix described in
Section~\ref{sec:impedance_Schur} is much more
suitable to use in practice within the framework of the
proposed method together with cross-points.
The use of the more involved communication matrix \(\Pi\)
computed by solving the projection problem is therefore worthwhile considering
in practice.

\paragraph{Convergence history of iterative algorithms}

\begin{figure}[htp]
    \centering
    \begin{subfigure}{0.49\textwidth}
        \centering
        \includegraphics[width=\textwidth]{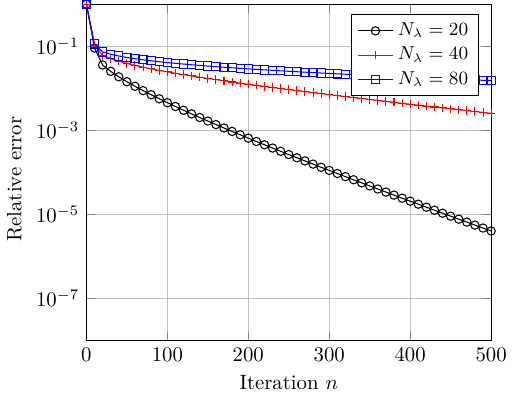}%
        \caption{Despr\'es \(\bfT\) matrix (Richardson).}\label{fig:MH_jacres_6_d}
    \end{subfigure}
    \begin{subfigure}{0.49\textwidth}
        \centering
        \includegraphics[width=\textwidth]{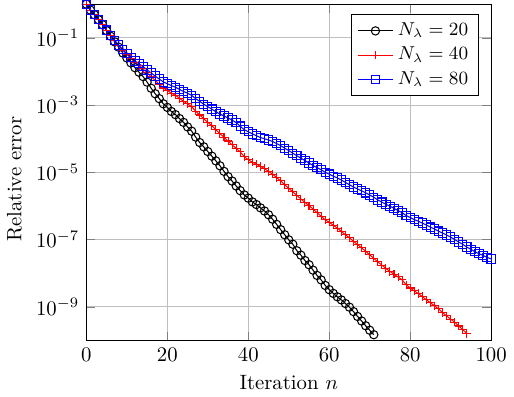}%
        \caption{Despr\'es \(\bfT\) matrix (\textsc{gmres}).}\label{fig:MH_gmresres_6_d}
    \end{subfigure}
    \begin{subfigure}{0.49\textwidth}
        \centering
        \includegraphics[width=\textwidth]{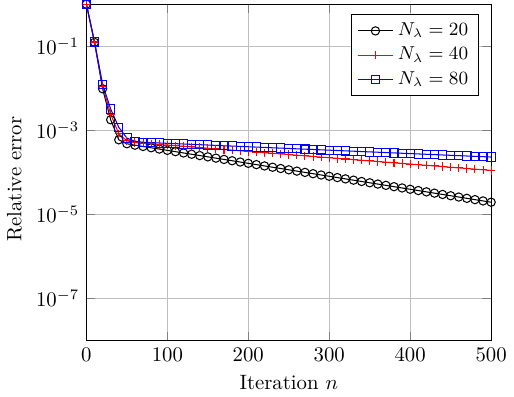}%
        \caption{Interface based non-local \(\bfT\) (Richardson).}\label{fig:MH_jacres_6_p}
    \end{subfigure}
    \begin{subfigure}{0.49\textwidth}
        \centering
        \includegraphics[width=\textwidth]{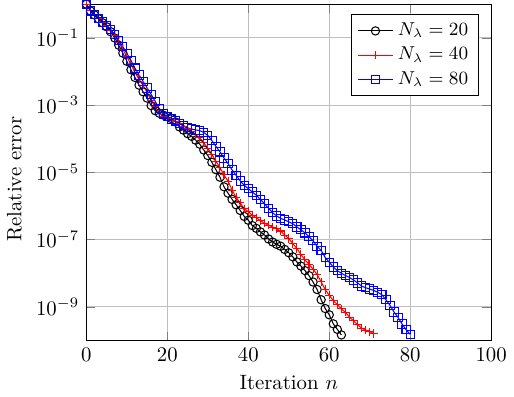}%
        \caption{Interface based non-local \(\bfT\) (\textsc{gmres}).}\label{fig:MH_gmresres_6_p}
    \end{subfigure}
    \begin{subfigure}{0.49\textwidth}
        \centering
        \includegraphics[width=\textwidth]{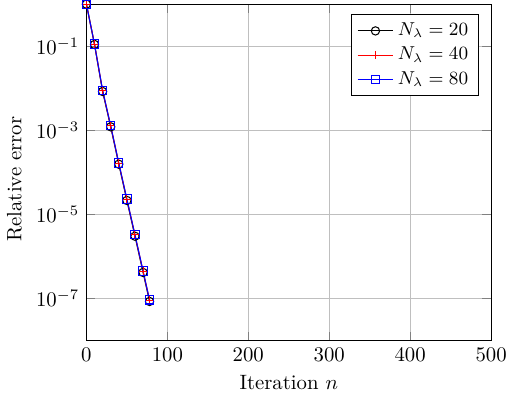}%
        \caption{Subdomain based non-local \(\bfT\) (Richardson).}\label{fig:MH_jacres_6_x}
    \end{subfigure}
    \begin{subfigure}{0.49\textwidth}
        \centering
        \includegraphics[width=\textwidth]{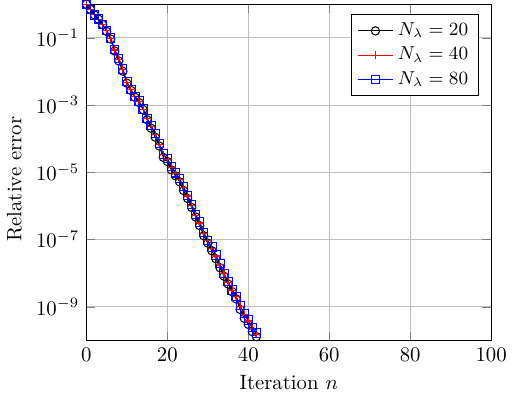}%
        \caption{Subdomain based non-local \(\bfT\) (\textsc{gmres}).}\label{fig:MH_gmresres_6_x}
    \end{subfigure}\caption{Convergence history for
    the Richardson algorithm (left) and \textsc{gmres}
    algorithm with a restart every 5 iterations (right).
    With \(\mJ=6\) subdomains and wavenumber \(\kappa=2\pi\).
    }\label{fig:MH_cv}
\end{figure}

We report in Figure~\ref{fig:MH_cv} the convergence histories
of the three domain decomposition methods
for the damped Richardson algorithm (left)
and for the \textsc{gmres} algorithm (right).
The results are provided for three different mesh refinements, indicated by
\(N_{\lambda}:=\lambda/h\) which is the number of points per wavelength
\(\lambda\) if \(h\) is the typical edge length.
The wavenumber is \(\kappa=2\pi\).
There are a total of respectively \(\#\Edge=4\;908\),
\(18\;180\) and \(71\;748\) degrees of freedom
for the three refinement considered 
\(N_{\lambda}=20\), \(40\) and \(80\).
We see the deterioration of the convergence of the iterative
algorithms with the mesh refinement when the Despr\'es transmission conditions
are used.
This is a common feature to transmission matrices based on local operators.
When the interface based non-local \(\bfT\) is used, we also see
a deterioration of the convergence with mesh refinement, albeit less pronounced.
Such observations were already reported in previous works~\cite{zbMATH07197799,
parolin:tel-03118712, claeys2020robust}.
On the contrary, the new approach based on a subdomain based non-local
\(\bfT\) that we developed exhibits a perfectly uniform convergence with
respect to the mesh size and converges faster than the other two strategies.

\paragraph{Eigenvalues of the iteration matrix}

\begin{figure}[htp]
    \centering
    \begin{subfigure}[b]{0.49\textwidth}
        \centering
        \includegraphics[width=\textwidth]{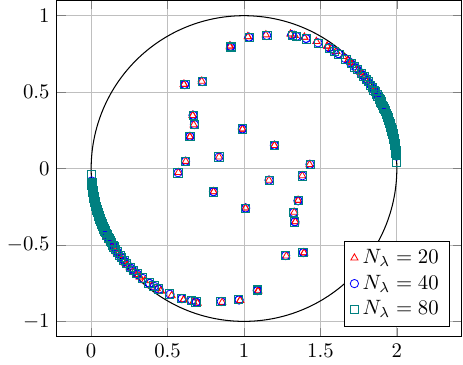}%
        \caption{Despr\'es \(\bfT\) matrix.}\label{fig:MH_vpA_6_d}
    \end{subfigure}
    \begin{subfigure}[b]{0.49\textwidth}
        \centering
        \includegraphics[width=\textwidth]{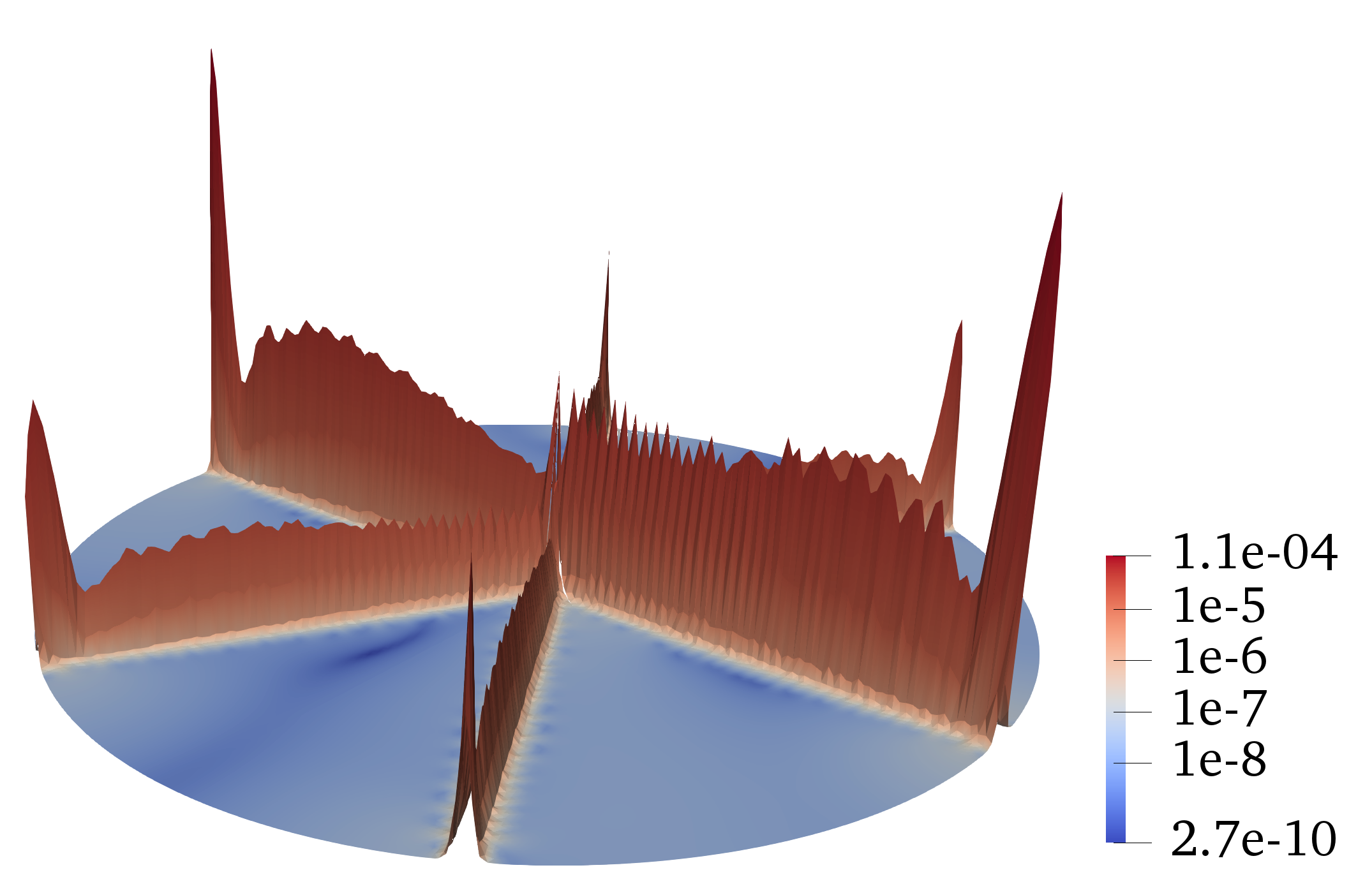}%
        \vspace{0.025\textheight}
        \caption{Despr\'es \(\bfT\) matrix.}\label{fig:MH_error_d_3D}
    \end{subfigure}
    \begin{subfigure}[b]{0.49\textwidth}
        \centering
        \includegraphics[width=\textwidth]{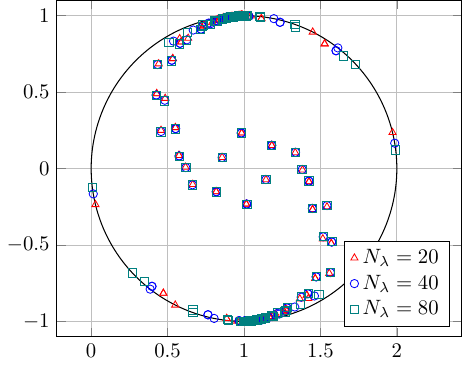}%
        \caption{Interface based non-local \(\bfT\).}\label{fig:MH_vpA_6_p}
    \end{subfigure}
    \begin{subfigure}[b]{0.49\textwidth}
        \centering
        \includegraphics[width=\textwidth]{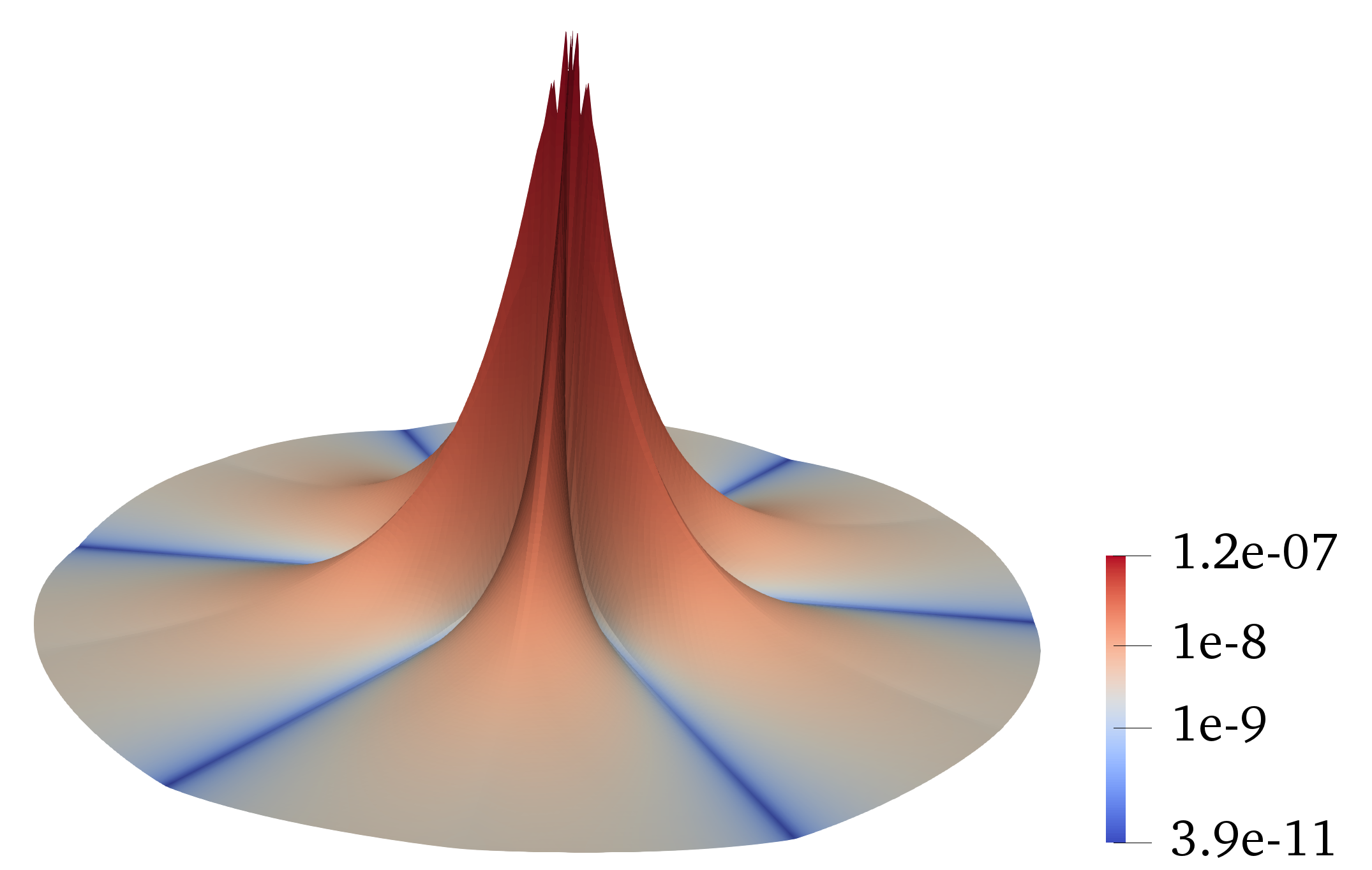}%
        \vspace{0.025\textheight}
        \caption{Interface based non-local \(\bfT\).}\label{fig:MH_error_p_3D}
    \end{subfigure}
    \begin{subfigure}[b]{0.49\textwidth}
        \centering
        \includegraphics[width=\textwidth]{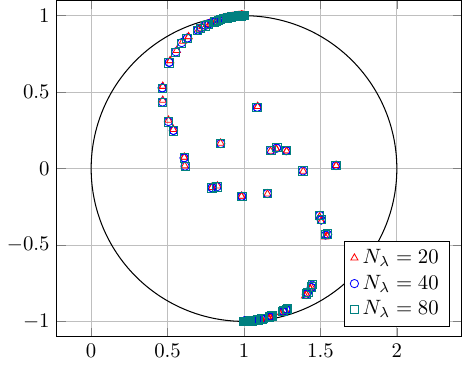}%
        \caption{Subdomain based non-local \(\bfT\).}\label{fig:MH_vpA_6_x}
    \end{subfigure}
    \begin{subfigure}[b]{0.49\textwidth}
        \centering
        \includegraphics[width=\textwidth]{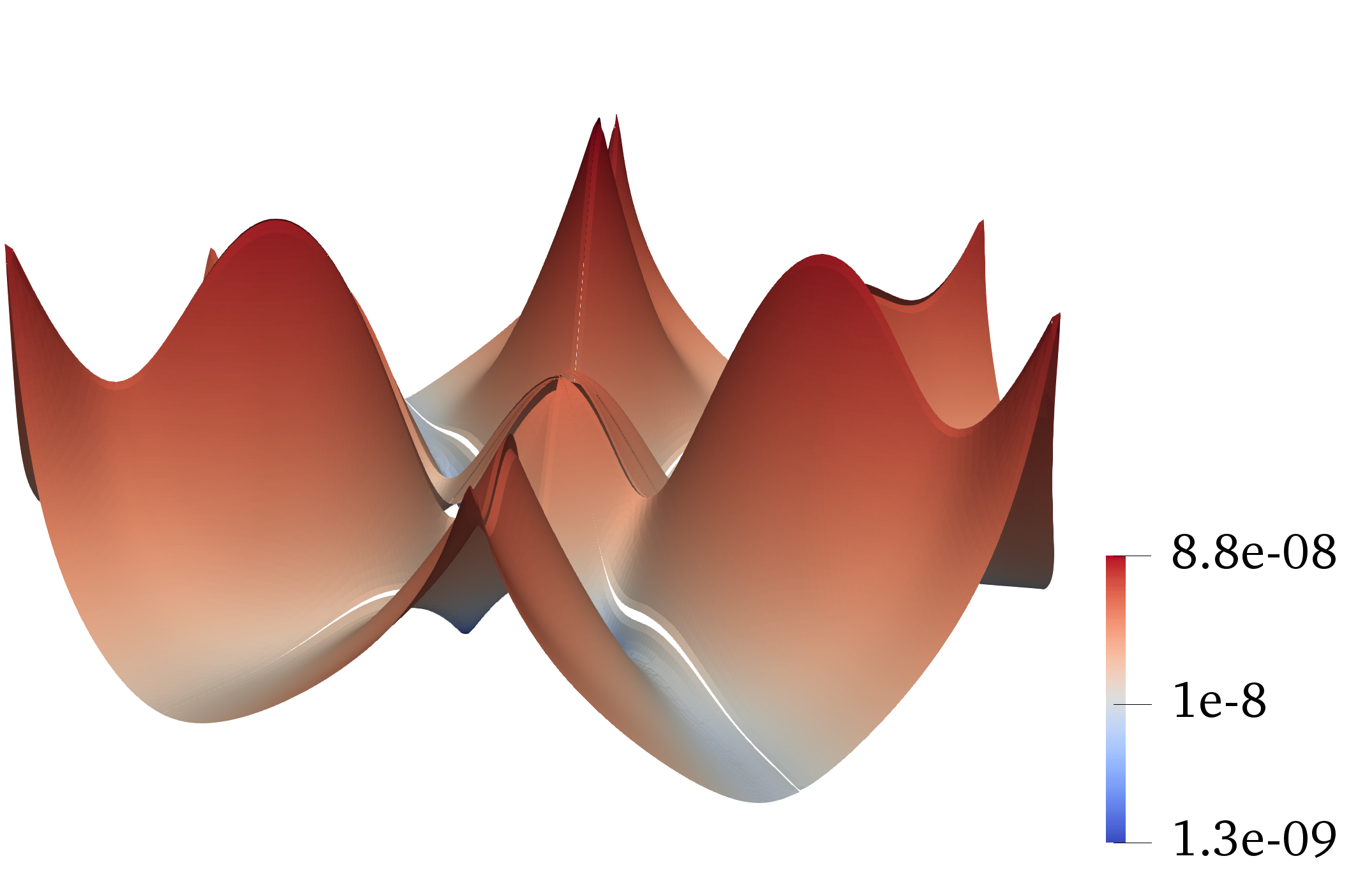}%
        \vspace{0.025\textheight}
        \caption{Subdomain based non-local \(\bfT\).}\label{fig:MH_error_x_3D}
    \end{subfigure}\caption{
    Eigenvalues of the iteration matrices $\bfId+\bfPi\bfS$ (left) 
    and nature of the error (right).
    The absolute value of the error on the solution is represented as the
    elevation (after linear interpolation on the nodes of the mesh).
    Different magnification factors are used for the three
    cases.
    With \(\mJ=6\) subdomains and wavenumber \(\kappa=2\pi\).
    }\label{fig:MH_vpa_error}
\end{figure}

To try to understand better those results, 
we report in Figure~\ref{fig:MH_vpa_error} (left)
the eigenvalues of the \emph{iteration matrices} $\bfId+\bfPi\bfS$ 
that are involved in the three domain decomposition methods.
When the Despr\'es transmission conditions are used, we observe an accumulation
close to the origin which will harm the convergence of both the \textsc{gmres}
and the damped Richardson algorithms.
When the interface based non-local \(\bfT\) is used, we see that
the clusters are near the two points \((1,1)\) and \((1,-1)\), which
demonstrates that the evanescent modes are well taken into account.
We see however a few isolated eigenvalues, close to the shifted unit
circle, which seem to get closer to the origin as the mesh is refined.
In contrast, with the subdomain based non-local \(\bfT\), the eigenvalues seem
to be uniformly bounded away from the critical points.

\paragraph{Nature of the error}

We represent in Figure~\ref{fig:MH_vpa_error} (right) the distribution of the error
between the exact discrete solution and the discrete solution (obtained with
the damped Richardson algorithm).  More precisely, the absolute value of the error
is represented as the elevation along the \(z\)-axis, after linear
interpolation on the nodes of the mesh.
For a better representation, the magnification factor is different for each
figure, as indicated by the actual maximum and minimal values of the error on
the colorbar.
The convergence is stopped before machine precision is reached.  In some sense,
the nature of the remaining error gives us insight on the components that are
troublesome for the convergence.

When Despr\'es transmission conditions are used, we see that the error
is highly concentrated along each interface and decreases very rapidly away from
them.
The most likely interpretation is that the main components in the error consist
in some sense of ``evanescent waves''.
Note also that the ratio between the maximum and minimum values of the error is
very large.

In contrast, the error is highly peaked at the cross-point and (slowly)
decreasing away from it when the interface based non-local \(\bfT\) is used.
The transmission interfaces seem less visible.
Note also that the ratio between the maximum and minimum values of the error is
much smaller than for the Despr\'es transmission conditions.

As for the subdomain based non-local \(\bfT\), the error is more
evenly distributed in the domain, albeit slightly accumulating near the
interfaces. 
More importantly, no accumulation of the error at the cross-point can be
observed in contrast to the result using also a non-local \(\bfT\)
but with the standard exchange matrix (Figure~\ref{fig:MH_error_p_3D}).

\paragraph{Influence of the choice of \(\Omega'\)}

\begin{figure}[ht]
    \centering
    \begin{subfigure}{0.4\textwidth}
        \centering
        \includegraphics[width=\textwidth]{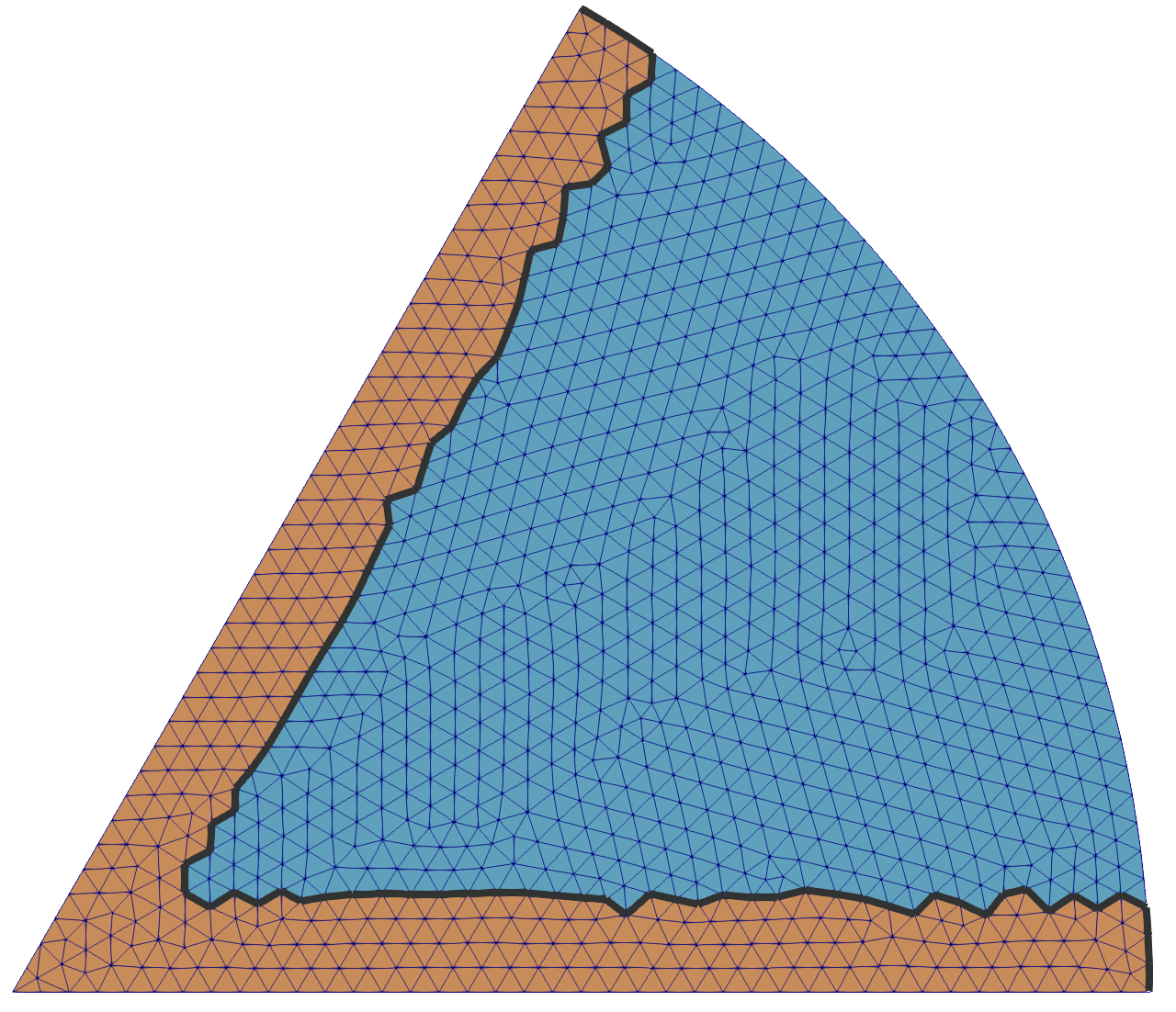}%
    \caption{\(N_{\lambda}=40\).}\label{fig:meshT_0040_0006}
    \end{subfigure}
    \hfill
    \begin{subfigure}{0.4\textwidth}
        \centering
        \includegraphics[width=\textwidth]{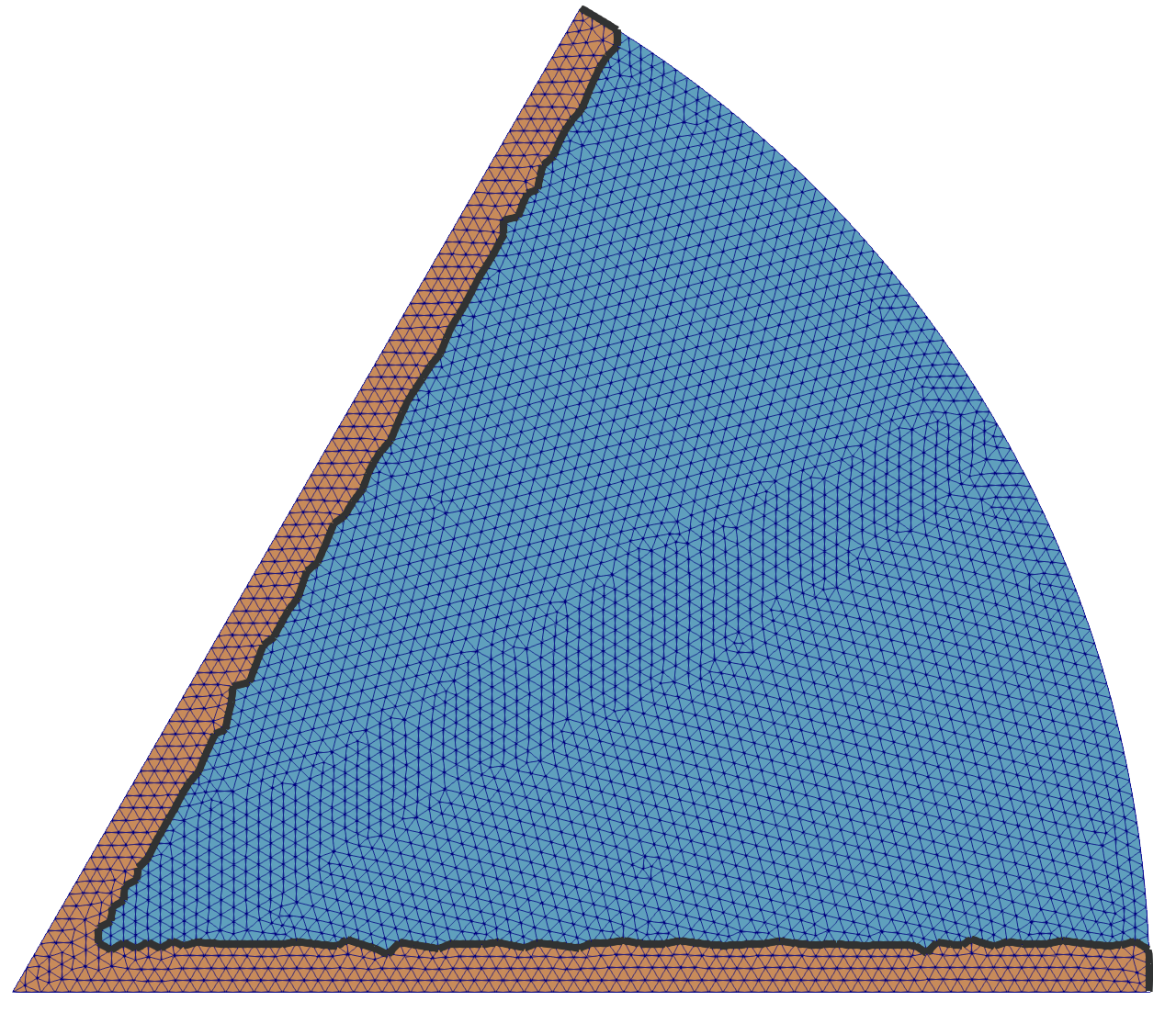}%
    \caption{\(N_{\lambda}=80\).}\label{fig:meshT_0080_0006}
    \end{subfigure}\caption{Definition of \(\Omega_{j}'\) (light brown region)
    used in the test case of Figure~\ref{fig:influence_omega_prime}.
    \(\Omega_{j}\) is the union of the two (light brown and light blue) regions.
    The dark edges consist of the domain of the boundary term
    in~\eqref{eq:schur_complement_bilinear_form}.
    Here \(J=6\), the other 5 subdomains can be obtained by rotation.}\label{fig:influence_omega_prime_meshes}
\end{figure}

\begin{figure}[ht]
    \centering
    \begin{subfigure}{0.49\textwidth}
        \centering
        \includegraphics[width=\textwidth]{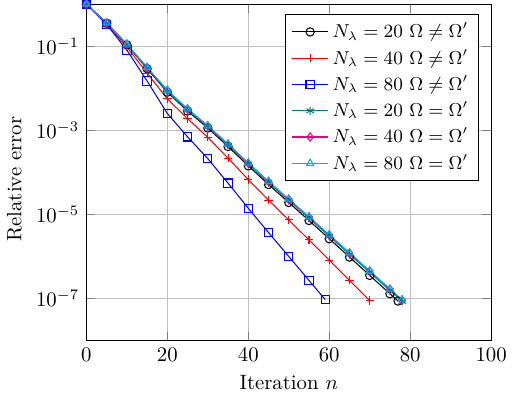}%
    \caption{Richardson algorithm.}\label{fig:jacres_0006_x_cvplot_Tcreux}
    \end{subfigure}
    \begin{subfigure}{0.49\textwidth}
        \centering
        \includegraphics[width=\textwidth]{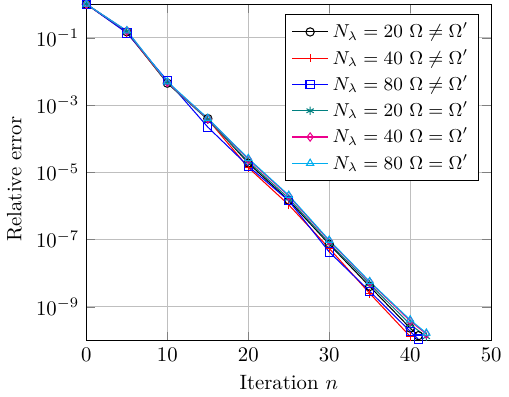}%
    \caption{\textsc{gmres} algorithm (restart 5).}\label{fig:gmresres_0006_x_cvplot_Tcreux}
    \end{subfigure}\caption{Influence of the choice of \(\Omega'\)
    on the convergence.
    Wavenumber \(\kappa=2\pi\).}\label{fig:influence_omega_prime}
\end{figure}

We finally investigate for this particular test case the influence of the choice
of the domain \(\Omega'\) that intervenes in the definition of the Schur
complement based matrix \(\bfT\), as defined in
Section~\ref{sec:impedance_Schur}.
The domain of the auxiliary problem \(\Omega_{j}'\) is represented in
Figure~\ref{fig:influence_omega_prime_meshes}.
It consists of the mesh cells that are within a distance of \(3h_{\max}\)
from the transmission boundary, where \(h_{\max}\) is the maximum edge length in
the triangulation.
The convergence results are given in
Figure~\ref{fig:influence_omega_prime}.
We see that using much smaller domains, concentrated in the vicinity of the
transmission boundaries has a minor (yet positive for the Richardson algorithm
in this particular case) effect on the convergence.
The computational cost of the matrix assembly is however greatly reduced.
This can be explained from the fact that we solve elliptic problems with 
a source term defined on the transmission boundary.
The solution is then mainly concentrated in the vicinity of this boundary.
Note that the boundary term in~\eqref{eq:schur_complement_bilinear_form}
is empirically found to be a crucial ingredient to obtain this result.
A modal analysis in a simple geometry as well as additional numerical
experiments regarding the choice of \(\Omega'\) can be found
in~\cite[Chap.~8]{parolin:tel-03118712}.

\subsection{Stability}

\begin{figure}[htp]
    \centering
    \includegraphics[width=0.45\textwidth]{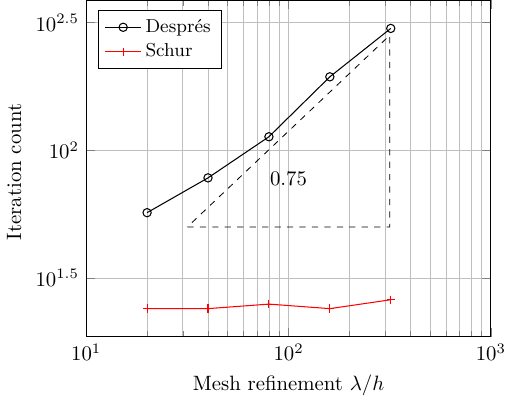}%
    \includegraphics[width=0.45\textwidth]{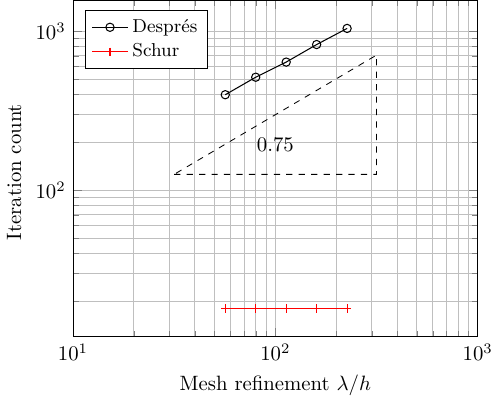}%
    \caption{Number of \textsc{gmres} iterations (restart 20) with respect to mesh refinement
    \(\lambda/h\), for 2D (left) and 3D (right)
    configurations.
    Wavenumber \(\kappa=1\).}\label{fig:stability}
\end{figure}

We investigate now further the robustness of the proposed approach with respect
to the mesh discretization, in particular with respect to mesh refinement now
both in 2D and 3D.
The refinement, namely decreasing the typical edge length \(h\), is uniform in
the domain \(\Omega\).
In the remainder of this manuscript and in contrast to the previous experiment,
the domain \(\Omega\) (a disk in 2D and a ball in 3D) will be partitioned using
an automatic graph partitioner. 
For the following results, there are \(\mJ=4\) subdomains in 2D and \(\mJ=32\)
subdomains in 3D and the wavenumber is \(\kappa=1\).
There are a total of \(\#\Edge=113\;627\) degrees of freedom in 2D and 
\(\#\Edge=137\;899\) degrees of freedom in 3D for the finest refinement.
The results are reported in Figure~\ref{fig:stability}.

We observe a quasi-linear increase in the number of iterations required to
reach a set tolerance for the Despr\'es \(\bfT\) matrix.
This is in stark contrast with the results using the Schur complement approach
which are completely immune to the mesh refinement.
Such an effect, which was already observed in the previous experiment, is
expected and not new,
see~\cite{zbMATH07197799,zbMATH01488462,parolin:tel-03118712}
and in particular the numerical analysis and numerical experiments
of~\cite{claeys2020robust} obtained in the acoustic setting.
In fact, it is one of the core strength of the approach based on the use of
underlying non-local operators in transmission conditions.

We shall point out however that in previous
works~\cite{zbMATH07197799,zbMATH01488462}
such an effect was observed, and as a matter of fact rigorously proved, only
in absence of cross-points in the partition.
Notice that in this 3D configuration (as a matter of fact, in all 3D tests cases
considered in this paper) there are indeed cross-points, namely degrees of
freedom with multiplicity strictly larger than two,
i.e.\ attached to edges that are shared by at least three sud-domains (such
points form the so-called \emph{wire-basket}).
This feature, namely the robustness with respect to mesh refinement, even in
presence of cross-points, is precisely enabled by the somewhat unusual choice
of communication matrix based on the global projection that was described in
the previous sections, see Section~\ref{sec:proj_in_practice}.

\begin{figure}[htp]
    \centering
    \includegraphics[width=0.45\textwidth]{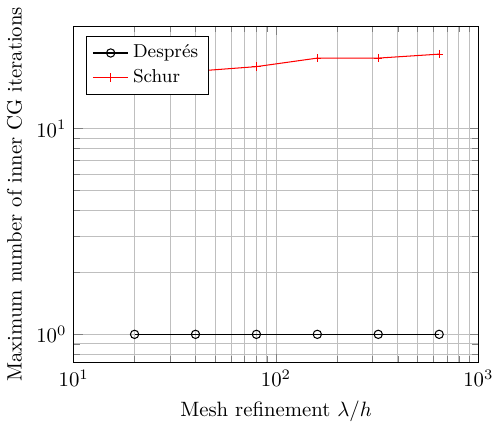}%
    \includegraphics[width=0.45\textwidth]{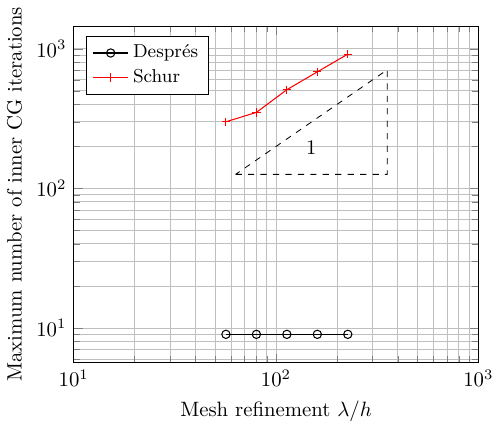}%
    \caption{Maximum number of iterations of the \emph{inner} preconditioned CG
    algorithm used to solve the projection problem with respect to mesh
    refinement \(\lambda/h\), for 2D (left) and 3D (right)
    configurations.}\label{fig:stability-cgmax}
\end{figure}

Besides, we report in Figure~\ref{fig:stability-cgmax} the number of iterations
of the \emph{inner} preconditioned CG algorithm that is used to solve the global
projection problem on the skeleton, see~\eqref{NeumannNeumannStrategy}.
We stress that these iteration counts do \emph{not} correspond to the
\emph{outer} iterations of the \textsc{gmres} algorithm that is still used to solve the
skeleton problem~\eqref{SkeletonFormulation}.

In 2D, we observe that a moderate number of iterations is required to solve the
projection problem using the transmission matrix based on a Schur complement.
It is moreover stable with mesh refinement.
For the Despr\'es \(\bfT\) matrix, we report exactly one iteration regardless of the
mesh refinement.
This is due to the fact that there are no degrees of freedom attached to the
cross-points in the two-dimensional configuration.
As a result, the linear system involved in the projection problem is actually
diagonal and there is no need to use the PCG algorithm in this particular case.

In 3D, we observe that a moderate number of iterations, stable with mesh
refinement, is required to solve the projection problem in the Despr\'es case.
This is expected since now there are actually degrees of freedom on the
junctions lines (or wire-basket) shared by at least three subdomains.
In contrast, we observe a linear growth of the number of iterations for the
Schur complement based approach.
Such a strong effect was not observed in the acoustic setting,
see~\cite{parolin:tel-03118712}.
It turns out that for the Maxwell setting, a more involved (auxiliary space)
preconditioning approach, based on a suitable Helmholtz-type splitting of edge
element vector fields, is necessary to tackle this issue~\cite{Hiptmair2006,
Hiptmair2007} but was not further explored in this first work.

\subsection{Influence of the number of subdomains}

\begin{figure}[htp]
    \centering
    \includegraphics[width=0.45\textwidth]{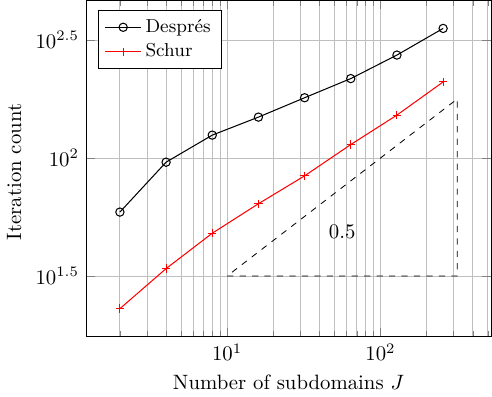}%
    \includegraphics[width=0.45\textwidth]{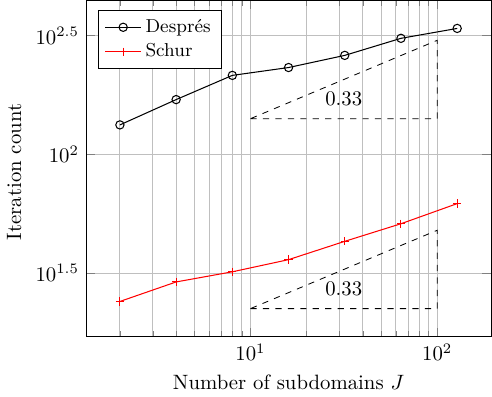}%
    \caption{Number of \textsc{gmres} iterations (restart 20) with respect to the number of
    subdomains \(\mJ\), for 2D (left) and 3D (right)
    configurations.
    Wavenumber \(\kappa=1\).}\label{fig:weak_scaling}
\end{figure}

We study now for both the 2D and 3D configurations the influence of the number of
subdomains \(\mJ\) on the number of iterations to reach a set tolerance with a
domain \(\Omega\) growing in size.
Specifically, the size of the domain is chosen to grow like \(\mJ^{1/d}\) where
\(d\) is the dimension of the ambient space, in order to keep a fixed size (in
terms of the number of degrees of freedom) for the local subproblems.
In 2D the domain is a disk of radius increasing from \(R=\sqrt{2}\) to
\(R=16\) as the number of subdomains increases from \(\mJ=2\) to \(\mJ=256\). 
In 3D the domain is a sphere of radius increasing from \(R=1\) to
\(R=4\) as the number of subdomains increases from \(\mJ=2\) to \(\mJ=128\).
In both cases, the wavenumber is \(\kappa=1\).
Notice that for this test case, despite the fact that the size of the problem
increases, the number of points per wavelength is kept constant. As a result
the pollution effect is not taken into account here.
There are a total of \(\#\Edge=113\;627\) degrees of freedom in 2D and 
\(\#\Edge=49\;877\) degrees of freedom in 3D for the largest \(\mJ\).
The results are provided in Figure~\ref{fig:weak_scaling}.

The growth of the number of iteration to reach the set tolerance also appears
to scale like \(\mJ^{1/d}\) and the phenomenon seems to apply to all the
transmission matrices considered.
This non-optimality is expected and can be understood in this wave propagation
context from the fact that the waves (hence the information) need to travel
longer distances as the size of the global domain increases.
Such an observation motivates the search for optimal solvers immune to this effect,
for instance using multi-level techniques and coarse spaces somehow mimicking
algorithms used for elliptic systems.
However, in this work, we did not pursue in this direction.

\subsection{Influence of the frequency}

\begin{figure}[htp]
    \centering
    \includegraphics[width=0.45\textwidth]{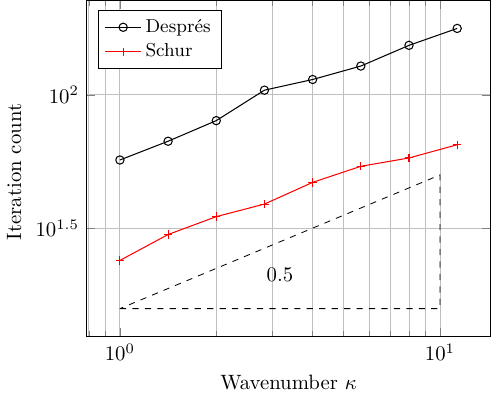}%
    \includegraphics[width=0.45\textwidth]{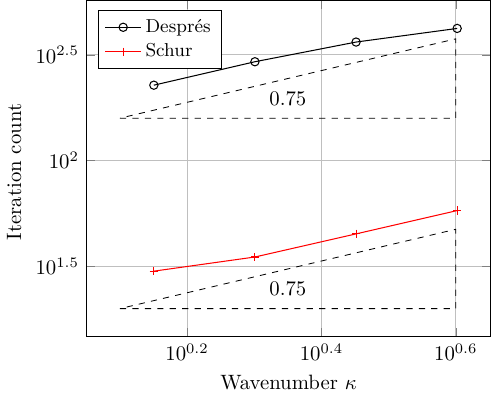}%
    \caption{Number of \textsc{gmres} iterations (restart 20) with respect to the wavenumber
    \(\kappa\), for 2D (left) and 3D (right)
    configurations.}\label{fig:pollution}
\end{figure}

We now study the dependency of the iteration counts with respect to the
wavenumber \(\kappa\).
To take the pollution effect into account, the mesh is refined as the frequency
increases. Since we are using low order finite elements,
we need to keep the quantity \(\kappa^{3}h^{2}\) fixed throughout the
computations to counter the pollution effect.
Here \(h\) denotes the typical edge length in the mesh.
In both the 2D and 3D configurations, this quantity is fixed to
\((2\pi)^2/400\) in order to have at least \(20\) points per wavelength for the
smallest wavenumber considered.
The domain \(\Omega\) is partitioned into \(\mJ=4\) subdomains in 2D and
\(\mJ=16\) subdomains in 3D.
There are a total of \(\#\Edge=160\;947\) degrees of freedom in 2D and 
\(\#\Edge=374\;889\) degrees of freedom in 3D for the largest wavenumber.
The results are reported in Figure~\ref{fig:pollution}.

As the wavenumber \(\kappa\) increases, the discrete (as well as the
continuous) problem gets harder (the condition number of the original
undecomposed matrix increases).
We notice a sub-linear increase of the number of iterations to reach the set
tolerance for all transmission matrices studied.
The increase seems to be stronger in the 3D configuration.

\subsection{Domain heterogeneity}\label{sec:num_heterogeneity}

To conclude this section on numerical experiments we present a more involved
test case with more complicated medium of propagation.
The objective is to illustrate the robustness of the proposed approach.
Specifically we consider three types of propagative medium
in our usual unit disk in 2D and unit ball in 3D.

\begin{figure}[htp]
    \centering
    \includegraphics[height=0.2\textheight]{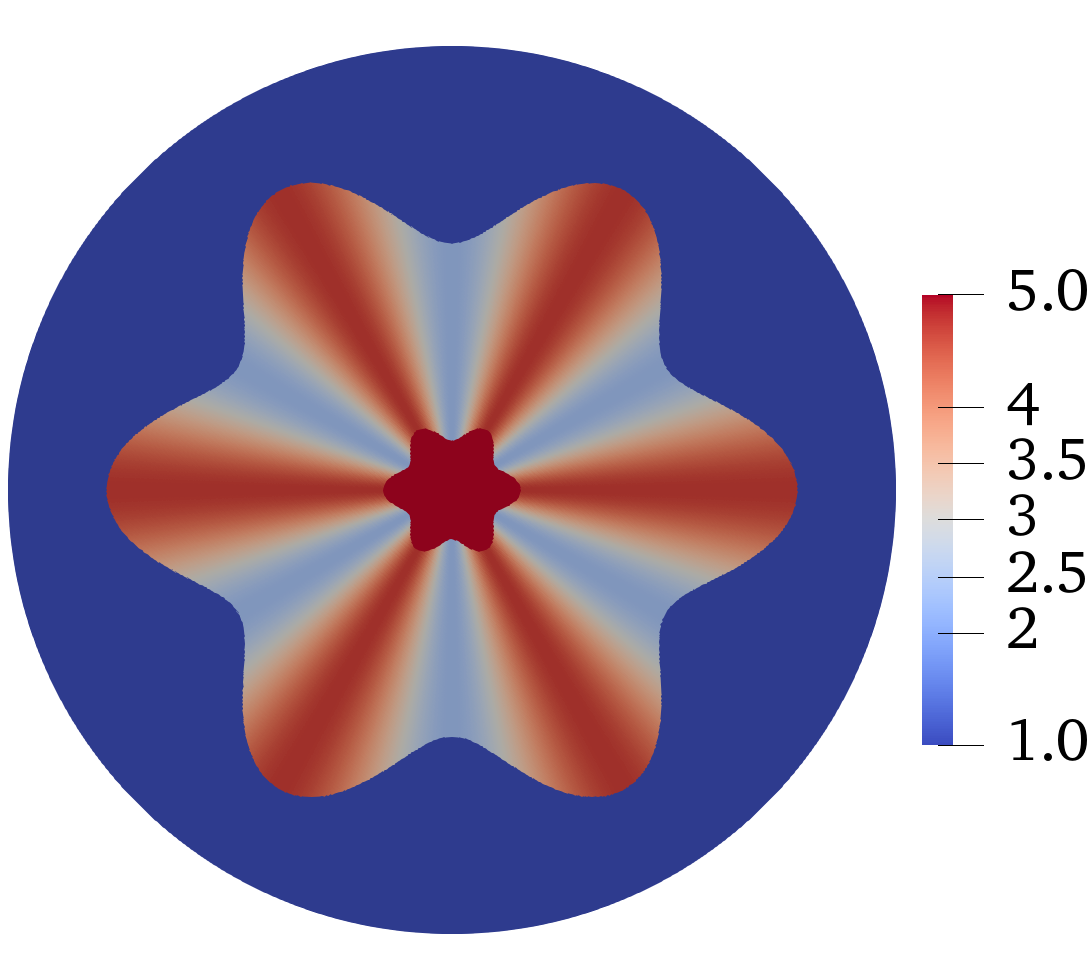}%
    \includegraphics[height=0.2\textheight]{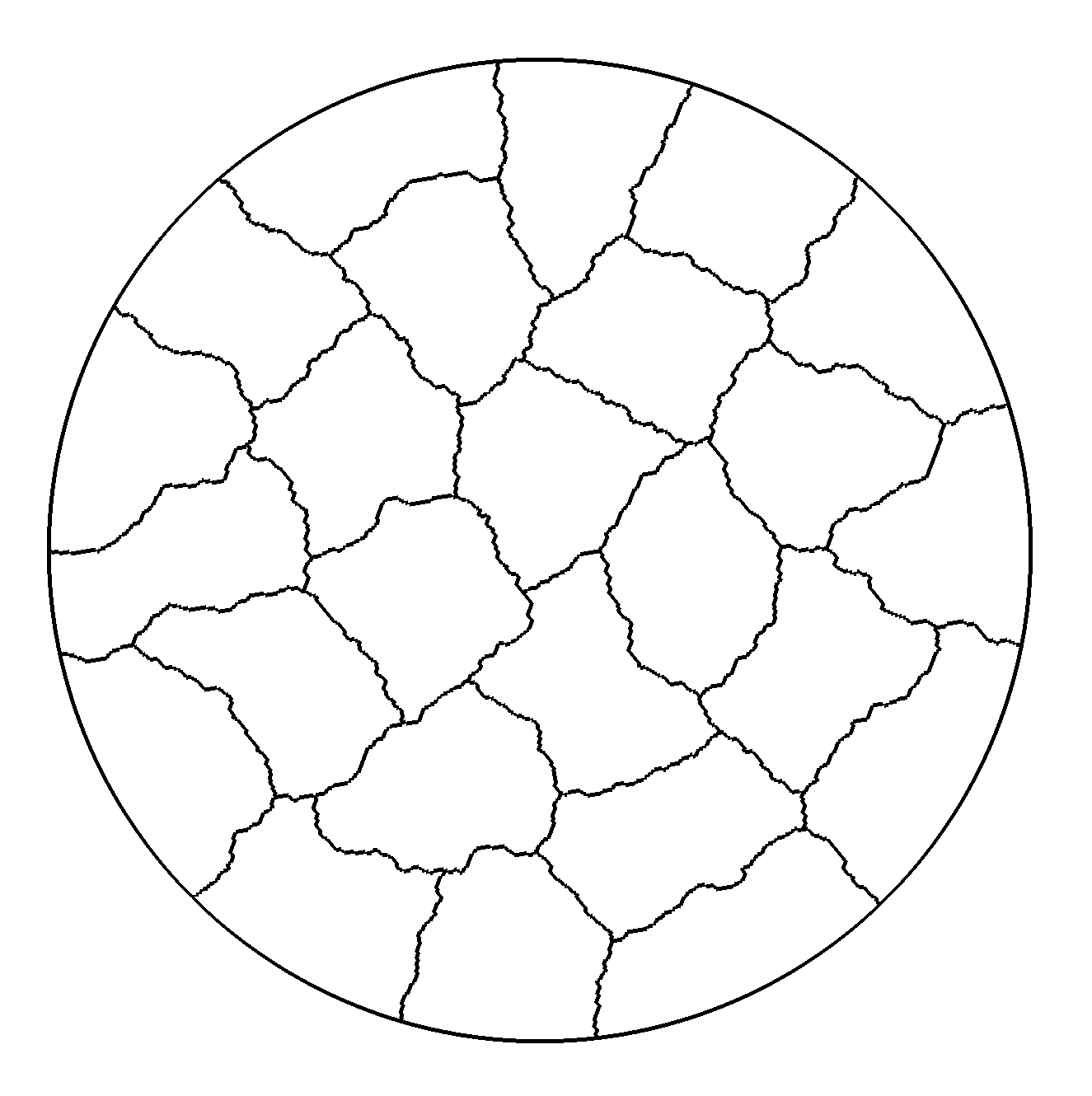}
    \includegraphics[height=0.2\textheight]{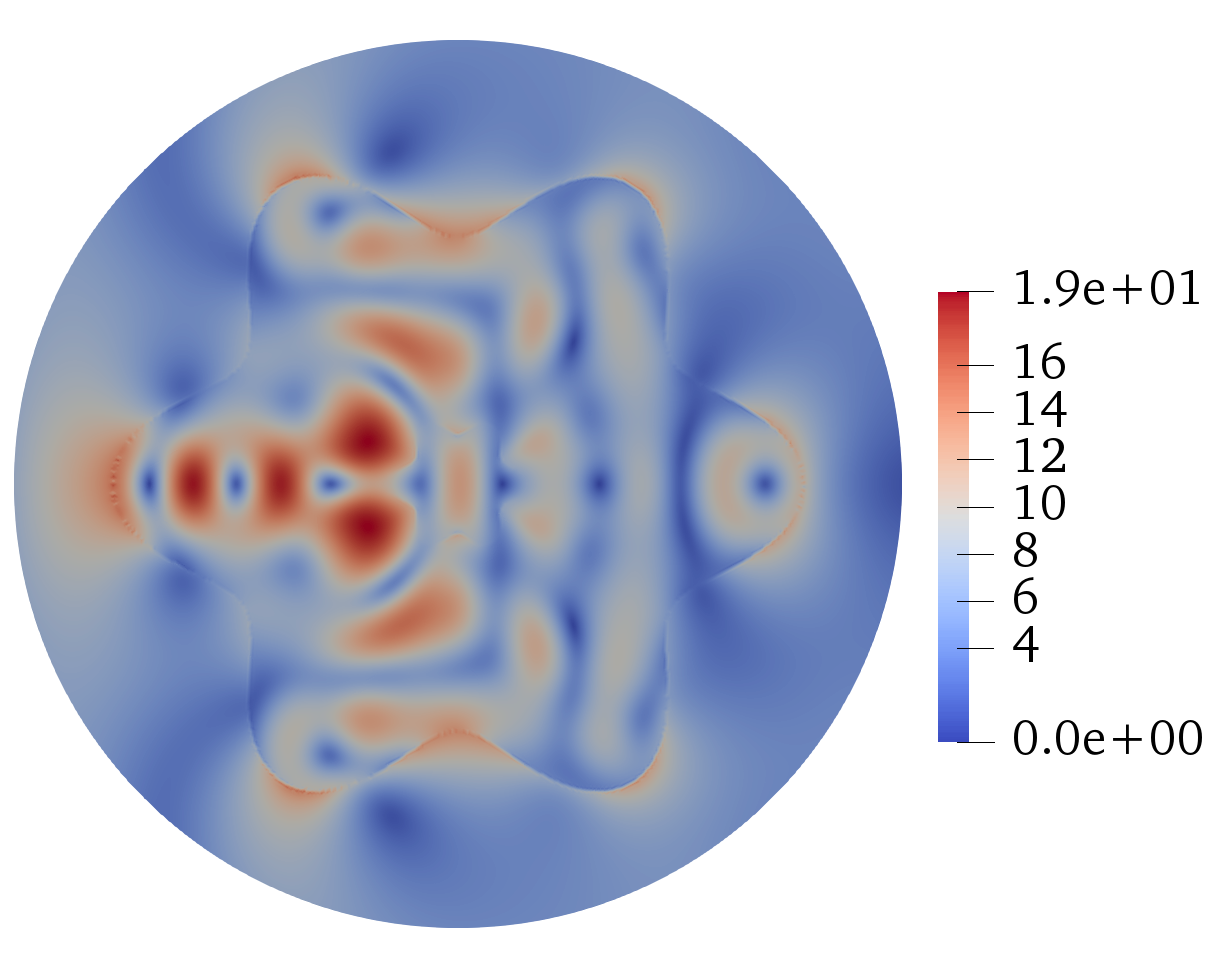}%
    \caption{Heterogeneous medium profile for the coefficient
    \(\check{\mu}_{r}\) (top left), skeleton of the partition (top
    right) and modulus of the solution for the purely propagative
    heterogeneous medium (bottom).}\label{fig:mu-profile_and_solution}
\end{figure}

The first medium is heterogeneous and purely propagative.
If \((r,\theta) \in [0,\, +\infty) \times [0,\, 2\pi)\)
and \((r,\varphi,\theta) \in [0,\, +\infty) \times [0,\, \pi) \times [0,\, 2\pi)\)
are respectively the cylindrical
and spherical coordinates, the coefficients \(\mu_{r} = \check\mu_{r}\) and
\(\epsilon_r = \check\epsilon_r\) are defined as follows
\begin{equation*}
  \check\mu_r := 
  \begin{cases}
    2\Delta\mu, & r \leq \rho(\theta)/5,\\
    1+\Delta\mu\,\psi(\theta), & \rho(\theta)/5 < r \leq \rho(\theta),\\
    1, & \rho(\theta) < r,\\
  \end{cases}
  \qquad
  \check\epsilon_r := 
  \begin{cases}
    2\Delta\epsilon, & r \leq \rho(\theta)/5,\\
    1+\Delta\epsilon\,\psi(\theta), & \rho(\theta)/5 < r \leq \rho(\theta),\\
    1, & \rho(\theta) < r,\\
  \end{cases}
\end{equation*}
where we set \(\Delta\mu=5/2\), \(\Delta\epsilon=3/2\) and
\begin{equation*}
  \rho(\theta) := 1+\cos(6\theta)/2,
  \qquad
  \psi(\theta) := 2(1+\cos(6\theta)/6)/3,
  \qquad\forall \theta\in[0,\, 2\pi).
\end{equation*}
See the left panel of
Figure~\ref{fig:mu-profile_and_solution} for a representation of the profile of
the \(\check\mu_{r}\) coefficient in 2D.
The coefficients are therefore both varying inside the domain and have surface
discontinuities.
The relative impedance is set to \(\eta_{r}=1\).
The wavenumber is set to \(\kappa=5\) in 2D and \(\kappa=1\) in 3D.

The second medium is homogeneous, constructed by averaging the coefficients of
the previous medium.
Specifically we used \(\mu_r=\epsilon_r=1\) and the wavenumber is set to
\(\kappa=5\check{\kappa}\) in 2D and \(\kappa=\check{\kappa}\) in 3D where \(\kappa_r\) is
the product of the averages on the domain \(\Omega\) of \(\check\mu_r\) and
\(\check\epsilon_r\) defined previously.

Finally, the third medium considered is heterogeneous and dissipative,
constructed by adding a strictly positive imaginary part to the coefficients of
the propagative heterogeneous medium previously defined. Specifically we used
\(\mu_r = \check\mu_r (1+\imath/4)\)
and \(\epsilon_r = \check\epsilon_r (1+\imath/6)\).
The wavenumber is set to \(\kappa=5\) in 2D and \(\kappa=1\) in 3D.

To simplify the comparison and discussion we used the same mesh (and partition)
in the three cases. Despite the possible heterogeneity of the medium, the mesh is
uniform, constructed such that the typical edge length parameter is
\(h=\lambda/50\) (resp. \(h=\lambda/30\)) with \(\lambda=2\pi/(5\kappa_{r})\)
(resp. \(\lambda=2\pi/\kappa_{r}\)) in 2D (resp. 3D).
The domain \(\Omega\) is partitioned into \(\mJ=25\) subdomains in 2D and
\(\mJ=50\) subdomains in 3D, see Figure~\ref{fig:mu-profile_and_solution}.
Since we are using an automatic graph partitioner independently of the
definition of the medium under consideration, some interfaces between two
subdomains are cut by the surface discontinuities of the coefficients (in the
heterogeneous case).
There are a total of \(\#\Edge=432\;103\) degrees of freedom in 2D and 
\(\#\Edge=310\;615\) degrees of freedom in 3D.

A solution is represented in the right panel of
Figure~\ref{fig:mu-profile_and_solution}, which corresponds to the 2D
propagative and heterogeneous medium configuration.
As before, the source comes from an impinging plane wave (coming from the left
in the Figure~\ref{fig:mu-profile_and_solution}).
Notice that, due to the heterogeneity, the modulus of the solution is rather
large in some part of the domain.

\begin{figure}[htp]
    \centering
    \includegraphics[width=0.45\textwidth]{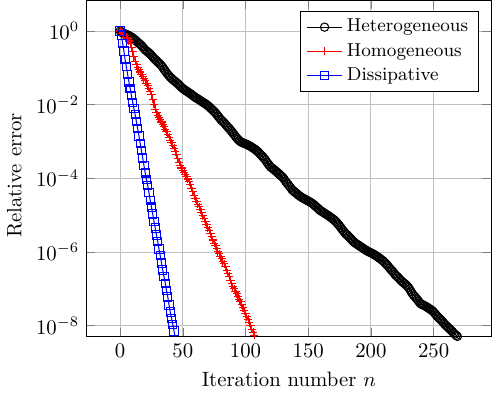}%
    \includegraphics[width=0.45\textwidth]{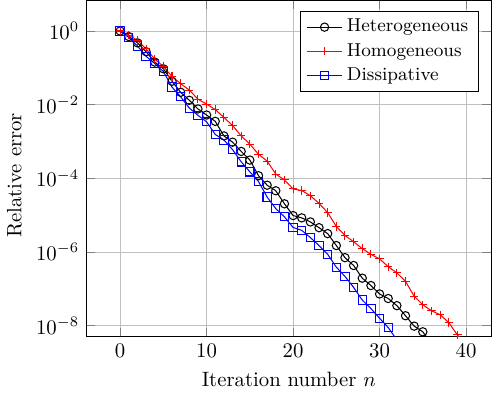}%
    \caption{Convergence history for heterogeneous, homogeneous and dissipative medium, for
    2D (left) and 3D (right) configurations. \textsc{gmres} algorithm (restart 20). }\label{fig:heterogeneity}
\end{figure}

We report in Figure~\ref{fig:heterogeneity} the convergence histories of the
\textsc{gmres} algorithm.
In the 2D case, we notice that a larger number of iterations is required in the
purely propagative heterogeneous medium (which is the notoriously more
difficult wave propagation problem) whereas the fastest convergence is achieved
in the dissipative scenario.
This is to be expected but we stress that the increase in the number of
iterations remains somewhat moderate.

In the 3D case, the convergence results are somewhat similar in the three
medium considered.
We explain this observation by noting that due to the relatively larger
frequency considered, the 2D test case corresponds to a more difficult wave
propagation problem than the 3D configuration.

\section*{Declarations}

\paragraph{Funding}
This work was supported by the project NonlocalDD funded by the French National
Research Agency, grant ANR--15--CE23--0017--01.

\paragraph{Conflict of interest}
The authors have no competing interests to declare.

\paragraph{Acknowledgments}
The authors would like to thank the two anonymous reviewers for their 
numerous relevant remarks which clearly improved the quality of this paper.

\bibliography{biblio}
\bibliographystyle{plain}

\end{document}